\documentclass[12pt]{iopart}

\usepackage{iopams}
\expandafter\let\csname equation*\endcsname\relax
\expandafter\let\csname endequation*\endcsname\relax
\usepackage{amsmath}
\usepackage{amsthm}
\usepackage{amssymb}
\usepackage{breqn}
\usepackage{tabu}
\usepackage{graphicx}
\usepackage[justification=centering]{caption}
\usepackage[colorinlistoftodos]{todonotes}

\newtheorem{thm}{Theorem}
\newtheorem{cor}[thm]{Corollary}

\newtheorem{defn}{Definition}
\newtheorem{prop}{Proposition}

\begin{document}

\title{Bounded Statistics}

\author{Pranava Chaitanya Jayanti}

\address{Department of Physics, University of Maryland, College Park, MD 20740, USA}
\ead{jayantip@umd.edu}

\author{Konstantina Trivisa}

\address{Department of Mathematics, University of Maryland, College Park, MD 20740, USA}
\ead{trivisa@math.umd.edu}

\vspace{10pt}
\begin{indented}
\item[]July 2018
\end{indented}

\begin{abstract}
If two probability density functions (PDFs) have values for their first $n$ moments which are quite close to each other (upper bounds of their differences are known), can it be expected that the PDFs themselves are very similar? Shown below is an algorithm to quantitatively estimate this “similarity” between the given PDFs, depending on how many moments one has information about. This method involves the concept of functions behaving “similarly” at certain “length scales”, which is also precisely defined. This technique could find use in data analysis, to compare a data set with a PDF or another data set, without having to fit a functional form to the data.
\end{abstract}

%
\vspace{2pc}
\noindent{\it Keywords}: moment problem, probability density functions, inverse problems, data analysis, characteristic functions \\
%
%
\maketitle
%
%

\section{Introduction} \label{introduction}

Modern scientific efforts often involve collecting, filtering, analysing and interpreting extremely large amounts of raw data. Famous examples from fundamental physics include the experiments at the Large Hadron Collider (LHC) and the Laser Interferometer Gravitational-Wave Observatory (LIGO). In the more applied sciences, one can think of complicated weather models that are tested for accuracy by comparison with meteorological data. Regression techniques are routinely used by engineers to fit very sophisticated equations of state to pressure-volume-temperature data of complex cryogenic mixtures. In all these cases, it would be of great use to develop a preliminary idea of the nature of a given data set or probability density function (PDF). This could help optimise computational costs.

To this end, one may ask: how much information regarding the underlying PDF can be extracted from the moments of a distribution? This is famously known as the moment problem \cite{Shohat1970TheMoments} and is a highly non-trivial question. It is a very important inverse problem with wide applications in data analysis. The question that is addressed in this article is a slight variant of the moment problem: is it possible to compare different distributions based on their moments alone, without having to explicitly construct their PDFs?

More precisely, the moment problem can be stated as follows: given a sequence of real numbers $\{\mu_k\}_{k=1,2,\dots}$, is it possible to find a PDF $f(x)$ such that $\{\mu_k\}$ are its moments, i.e., $\mu_k=\int_a^bx^k f(x) dx$? Depending on the domain of integration, the moment problem is classified into three categories: the Hausdorff problem which has a finite domain (usually taken to be $[0,1]$ with no loss in generality), the Stieltjes problem with a half-infinite domain $[0,\infty)$, and the Hamburger problem which spans the whole real line $(-\infty,\infty)$. These cases have been well-studied and conditions for existence and/or uniqueness of solutions have been formulated \cite{Shohat1970TheMoments}.

However, in many cases, an explicit construction of the PDF \cite{Talenti1987RecoveringMoments} may not be the goal, but it is done anyway. For instance, consider the ubiquitous situation of comparing two data sets. One way to do this is to fit PDFs using some regression techniques to the two data sets, and then compare the two PDFs. For the Hausdorff moment problem, this procedure of constructing the PDF from a finite set of moments is ill-defined in the sense of Hadamard, i.e., it is highly sensitive to the values of the moments due to the (exponentially) large condition number of the corresponding Hankel matrix \cite{Talenti1987RecoveringMoments}. It is also computationally very expensive, especially for large data sets. The computational cost can be decreased if an approximate form for the PDF can be deduced a priori, and this is the focus of this article. Given the first $n$ moments of each distribution/data set, or even just how far apart they are from each other, it is shown that the closeness (or similarity) of the underlying distributions can be commented upon in a very quantitative manner.

The moments of a distribution are directly related to the derivatives of the characteristic function of the PDF. Hence, the closeness of the moments of two distributions directly translates to the similarity of the derivatives (at $k=0$) of their characteristic functions. This fact is used, along with Taylor series expansions of the characteristic functions, to precisely characterize the “similarity” of the two distributions. In the case of a Hausdorff moment problem, uniqueness is guaranteed, i.e., \textit{if} there exists a function that solves the problem, it is unique \cite{Shohat1970TheMoments,Talenti1987RecoveringMoments}. In other words, for the case of a finite domain, the moments characterize the PDF completely. Thus, the two PDFs can be expected to be the same if and only if they have exactly the same moments, which is a trivial case. In the case that they are different, the PDFs may be expected to behave similarly only on certain “large” cut-off length scales, and differ significantly on the finer scales. This concept of “scale-wise similarity” is also precisely defined in the article. The result of this algorithm is the estimation of a cut-off length scale over which the two PDFs are similar (as defined by a user-specified tolerance), given how far apart the first few moments are.

The outline of this article is as follows: Section \ref{definitions} of this article defines some of the concepts used in the algorithm and Section \ref{estimates} establishes the estimates used in this work. Section \ref{existence & uniqueness} describes the algorithm, while Section \ref{infinite moments} consists of a brief discussion on the cases of bounded and unbounded support. Then, we illustrate the algorithm with some numerical examples in Section \ref{numerical examples}, following which the method is extended to higher dimensions in Section 7. Conclusions and scope for future work are presented in Sections 8 and 9, respectively.

\section{Definitions} \label{definitions}

\subsection{Distance metric} \label{distance metric}

Consider two absolutely integrable functions over a compact support, $f,g : [a,b]\rightarrow \mathbb{R}$. Define a distance metric between the two functions as:

\begin{equation}
        d(f,g) = \frac{\int_a^b |f(x)-g(x)| dx}{\int_a^b dx}
\end{equation}

\noindent Notice that this definition conforms to all the requirements of a distance metric.

\subsection{Characteristic function} \label{characteristic function}

The characteristic function (denoted by a capital letter) is the Fourier transform of a PDF. Here, and in what follows, the Fourier transform and its inverse are given by:

\begin{equation} \label{eq:characteristic function}
    F(k) = \int_a^b e^{ikx}f(x) dx \ , \ f(x) = \frac{1}{2\pi}\int_{-\infty}^{\infty} e^{-ikx}F(k) dk
\end{equation}

Assuming that the derivatives of the characteristic function exist, and that they all vanish at infinity, it can be easily shown that for a normalized PDF of a random variable $x$, the $n^{\text{th}}$ moment can be expressed as:

\begin{equation}
    \langle x^n \rangle = \int_a^b x^n f(x) dx = (-i)^n F^{(n)}(0)
\end{equation}

\noindent where $F^{(n)}(0) = \frac{d^n}{dk^n}F(k) \vert_{k=0}$.

Let the $n^{\text{th}}$ moments of two normalized PDFs $f(x)$ and $g(x)$ be denoted by $\langle x^n \rangle_f$ and $\langle x^n \rangle_g$, respectively. Define the bound on the difference of the $n^{\text{th}}$ moments  by:

\begin{equation}
    M_n \ge \left\vert \langle x^n \rangle_f - \langle x^n \rangle_g \right\vert = \left\vert F^{(n)}(0) - G^{(n)}(0) \right\vert
\end{equation}

\subsection{Low-pass filtering} \label{low-pass filtering}

In what follows, the notation in \cite{Frisch1995Turbulence:Kolmogorov} has been used to describe the large scale behaviour of functions. Filtering operators, one of which is extensively used in this article, are defined below.

\begin{defn}
    The action of the low-pass filtering operator $P_K$, associated with inverse length scale $K$, on an integrable function $f(x)$ is given by:
    
    \begin{equation} \label{eq:low-pass filter}
        P_K[f(x)] = f^<(x) = \frac{1}{2\pi}\int_{-K}^{K} e^{-ikx}F(k) dk
    \end{equation}
    
    \noindent where $F(k)$ is the Fourier transform of $f(x)$.
\end{defn}
\qed \\

The low-pass filtering operator has a “smoothing effect” in the sense that it removes any fluctuations of length scale less than $Κ^{-1}$. It is trivial to show that $P_K$ is a projection operator. The complement of the low-pass filter is called the high-pass filter, which removes the large scale behaviour and returns the small-scale features of the function.

\begin{defn}
    The action of the high-pass filtering operator $Q_K$, associated with inverse length scale $K$, on an integrable function $f(x)$ is given by:
    
    \begin{equation*}
        Q_K[f(x)] = (\mathbb{I} - P_K)[f(x)] = f^>(x) = \frac{1}{2\pi}\int_{|k|>K} e^{-ikx}F(k) dk
    \end{equation*}
\end{defn}
\qed \\

\begin{figure}
  \includegraphics[width=\linewidth]{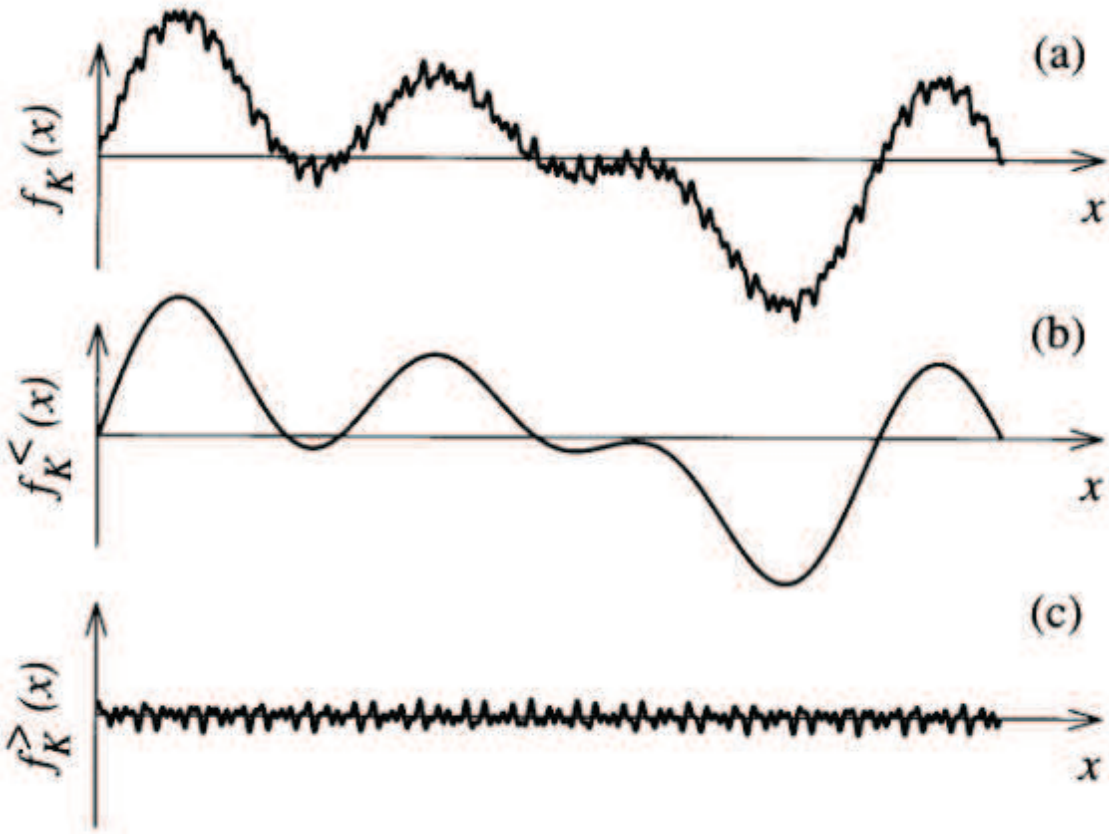}
  \caption{(a) Test function and action of (b) low-pass and (c) high-pass filters. (Reproduced with permission from Chapter 2 of \cite{Frisch1995Turbulence:Kolmogorov})}
  \label{fig:figure1}
\end{figure}

As an illustration, the action of these filters on a test function is shown in Figure \ref{fig:figure1}.

\section{Estimates of the distance} \label{estimates}

As stated in the introduction, we will use the characteristic relation between the moments and the derivatives of the characteristic functions to construct the required algorithm. Consider two PDFs $f(x)$ and $g(x)$ over the compact support $[0,1]$. Let their respective characteristic functions be $F(k)$ and $G(k)$. Assuming we have knowledge of the first $n$ moments (starting from the zeroth moment), we expand the characteristic functions in a Taylor series about $k=0$.

\begin{equation*}
    \begin{split}
        F(k) &= F^{(0)}(0) + \frac{F^{(1)}(0)}{1!}k + \frac{F^{(2)}(0)}{2!}k^2 + \dots + \frac{F^{(n-1)}(0)}{(n-1)!}k^{n-1} + \frac{F^{(n)}(t_Fk)}{n!}k^n\\
        G(k) &= G^{(0)}(0) + \frac{G^{(1)}(0)}{1!}k + \frac{G^{(2)}(0)}{2!}k^2 + \dots + \frac{G^{(n-1)}(0)}{(n-1)!}k^{n-1} + \frac{G^{(n)}(t_Gk)}{n!}k^n
    \end{split}
\end{equation*}

\noindent where $0\le t_F,t_G \le 1$. The last term in the expansion is the Lagrange form of the remainder in Taylor's theorem, which is applicable when the function being expanded as a series belongs to $C^{n-1}$ over the entire domain and is $n$-times differentiable. (For a rather straightforward proof that the characteristic function satisfies these conditions, see \ref{appendix_smoothness_characteristic function}).

Using these Taylor expansions, and the bounds on the difference of the moments (Section \ref{characteristic function}),

\begin{equation*}
\begin{split}
    |F(k) - G(k)| &\le \left\vert F^{(0)}(0) - G^{(0)}(0)\right\vert + \left\vert(F^{(1)}(0)-G^{(1)}(0))\frac{k}{1!}\right\vert + \left\vert(F^{(2)}(0)-G^{(2)}(0))\frac{k^2}{2!}\right\vert\\ 
    &+ \dots + \left\vert(F^{(n-1)}(0)-G^{(n-1)}(0))\frac{k^{n-1}}{(n-1)!}\right\vert + \left\vert(F^{(n)}(t_Fk)-G^{(n)}(t_Gk))\frac{k^n}{n!}\right\vert\\
\end{split}
\end{equation*}
    
\begin{equation}
    \begin{split}
    \Rightarrow |F(k) - G(k)| &\le M_0 + \frac{M_1}{1!}|k| + \frac{M_2}{2!}|k|^2 + \dots + \frac{M_{n-1}}{(n-1)!}|k|^{n-1}\\ 
    &\qquad \qquad \qquad + \frac{\left\vert F^{(n)}(t_Fk)-G^{(n)}(t_Gk)\right\vert}{n!}|k|^n\\
    \end{split}
\end{equation}

From \eqref{eq:low-pass filter}, we have the following inequality:

\begin{equation*}
    |f^<(x)-g^<(x)|=\left\vert \frac{1}{2\pi} \int_{-K}^{K} e^{-ikx}\left( F(k)-G(k) \right) dk\right\vert \le \frac{1}{2\pi}\int_{-K}^{K} |F(k)-G(k)|dk
\end{equation*}

By the definition of the distance metric and using this inequality,

\begin{equation*}
    \begin{split}
    d(f^<,g^<) &= \frac{\int_0^1 |f^<(x)-g^<(x)|dx}{\int_0^1 dx} \le \frac{\int_0^1 \frac{1}{2\pi}\int_{-K}^{K} |F(k)-G(k)|dk dx}{\int_0^1 dx}\\
    &\le \frac{1}{2\pi}\int_{-K}^{K} \Big\{ M_0 + \frac{M_1}{1!}|k| + \frac{M_2}{2!}|k|^2 + \dots + \frac{M_{n-1}}{(n-1)!}|k|^{n-1}\\ 
    &\qquad \qquad \qquad \qquad + \frac{\left\vert F^{(n)}(t_Fk)-G^{(n)}(t_Gk)\right\vert}{n!}|k|^n \Big\} dk \\
    \end{split}
\end{equation*}

Recalling \eqref{eq:characteristic function}, we have the following bound:

\begin{equation*}
    F^{(n)}(t_Fk) = i^n\int_0^1 x^n f(x) e^{it_fkx} dx \Rightarrow \left\vert F^{(n)}(t_Fk) \right\vert \le  \int_0^1 x^n f(x) dx = \langle x^n \rangle_f
\end{equation*}

Similarly, $\left\vert G^{(n)}(t_Gk) \right\vert \le \langle x^n \rangle_g$. Combining these results,

\begin{equation} \label{eq:preliminary distance}
    d(f^<,g^<) \le \frac{1}{\pi} \left\{ M_0K + \frac{M_1}{2!}K^2 + \dots + \frac{M_{n-1}}{n!}K^n + \frac{\langle x^n \rangle_f + \langle x^n \rangle_g}{(n+1)!}K^{n+1} \right\} 
\end{equation}

At this stage, all that remains is for us to estimate the unknown $\langle x^n \rangle_f$ (and $\langle x^n \rangle_g$) in terms of the lower moments, which are known. Depending on how much a priori information we have about the PDF $f(x)$, we may be able to estimate these unknown higher moments in various ways (see \ref{bounds on higher moments}). At this stage, we will consider the simplest case, where we know nothing about $f(x)$ except that it is a PDF, i.e., it is integrable. Since $0\le x\le 1$, using H\"older's inequality, it is easy to establish that

\begin{equation*}
    \langle x^n \rangle_f \le \langle x^m \rangle_f \ \forall \ 0\le m<n
\end{equation*}

Thus, the best (smallest) upper bound is when $m=n-1$. From \eqref{eq:preliminary distance} and this estimate,

\begin{equation} \label{eq:final distance}
    d(f^<,g^<) \le \frac{1}{\pi} \left\{ M_0K + \frac{M_1}{2!}K^2 + \dots + \frac{M_{n-1}}{n!}K^n + \frac{R_n}{(n+1)!}K^{n+1} \right\} 
\end{equation}

\noindent where $R_n = \langle x^{n-1} \rangle_f + \langle x^{n-1} \rangle_g$.\\

This means that the two PDFs have “similar” behaviour over length scales greater than 1/Κ, if the distance between the smoothed functions is less than some specified tolerance. The functions may, however, differ from each other significantly over small scales, which is expected since their moments are not exactly equal, but are only similar in the sense of their difference being bounded above.

\section{Existence and uniqueness of solutions} \label{existence & uniqueness}

Given a PDF, its closeness (at a certain length scale) to any other PDF can be estimated by specifying a tolerance for considering two functions to “behave similarly” and checking if there exists a length scale over which this is achieved. The question to be posed is the following:

\textit{For any $\varepsilon>0$, does there exist an inverse length scale $K$ such that $d(f^<,g^< )≤\varepsilon$, where $f^<$ and $g^<$ are the low-pass filtered functions as defined in \eqref{eq:low-pass filter}?}

This is answered by looking at the polynomial equation (in the variable $K$):

\begin{equation} \label{eq:polynomial in K}
    \frac{1}{\pi} \left\{ M_0K + \frac{M_1}{2!}K^2 + \frac{M_2}{3!}K^3 + \dots + \frac{M_{n-1}}{n!}K^n + \frac{R_n}{(n+1)!}K^{n+1}\right\} - \varepsilon = 0
\end{equation}

\noindent and calculating the inverse length scale $K$ which is a root of this equation. \eqref{eq:polynomial in K} is a polynomial equation of the $(n+1)^{\text{th}}$ degree, leading to the existence of $n+1$ solutions. Of these, the sought-after one is a real and positive value of $K$, preferably a unique root (to prevent further complications of having to choose between multiple solutions).

\begin{prop}
    There exists a unique positive solution to \eqref{eq:polynomial in K}, which characterizes the desired cut-off scale.
\end{prop}

\textit{Proof}: \eqref{eq:polynomial in K} is a polynomial equation of degree $n+1$ with all positive coefficients, except the constant term. Using Descartes’ rule of signs, it can be concluded that there is exactly one positive root, since exactly one sign change occurs throughout the polynomial. 
\qed \\

\section{Infinite/non-existent moments} \label{infinite moments}

\subsection{Bounded support} \label{bounded support}

For distributions defined on a bounded support, all moments exist regardless of the PDF. This can be shown as follows. (Since any bounded domain $[a,b]$ can be mapped to the compact set $[0,1]$, only the latter shall be considered, without any loss of generality.)

\begin{equation*}
    \langle x^n \rangle = \int_0^1 x^n f(x) dx \le \int_0^1 f(x) dx = 1
\end{equation*}

\noindent where the last equality follows from normalization, which is valid for all PDFs since they are Lebesgue-integrable (by definition). Thus, any PDF on a bounded support has all moments and they are all finite. Hence, the method described in this report is certainly applicable to PDFs with bounded support.

\subsection{Unbounded support} \label{unbounded support}

It is not necessary that all the moments, or even any moments, exist for a PDF with an unbounded support. For instance, while the normal distribution possesses all moments, the Cauchy distribution has none. If the $n^{\text{th}}$ moment of a distribution does not exist, it means that the $n^{\text{th}}$ derivative of the characteristic function does not exist as a limit at $k=0$. Similarly, the $n^{\text{th}}$ derivative of the characteristic function could blow up at $k=0$, signifying the blow up of the $n^{\text{th}}$ moment. In the latter case, the limit exists and is equal to $\infty$. \\

\noindent \textit{Note: However, it can be easily shown that for a PDF with an unbounded support, the existence of a finite $n^{\text{th}}$ moment implies the existence and finiteness of all moments that are lesser than $n$. For instance, $\forall \  0\le m<n$,}

\begin{align*}
    \langle x^n \rangle = \int_0^{\infty} x^n f(x) dx &= \int_0^1 x^n f(x) dx + \int_1^{\infty} x^n f(x) dx\\
    &\ge \int_0^1 x^n f(x) dx + \int_1^{\infty} x^m f(x) dx\\
    &= \int_0^1 (x^n - x^m) f(x) dx + \langle x^m \rangle
\end{align*}

\noindent\textit{The first term in the last step is clearly finite.}
\qed \\

If the higher derivatives of the characteristic functions become arbitrarily large, this corresponds to very large fine-scale fluctuations (this interpretation follows from the relation between the moments and the derivatives of the characteristic function). These fluctuations are not physical in the sense that they will prevent any kind of macroscopic equilibration from occurring, and hence characterize far-from-equilibrium processes in their transient states\footnote{Note that while turbulent flow is also (strictly speaking) far-from-equilibrium, there is a difference between the stationary and transient states of turbulent flow. In the former, there is statistical equilibrium, and the Kolmogorov spectrum is evidence of fluctuations decreasing with reducing length scales.}. Such states of these processes may not even have a well-defined PDF to begin with, and the method of this article is irrelevant to such extreme cases.

\section{Numerical examples} \label{numerical examples}

In this section, the above method is applied to some PDFs with compact support $[0,1]$. Two cases are considered: a “structured” PDF (with a closed-form expression) and an “unstructured" PDF (generated with random numbers). In both cases, various realisations of the PDFs are compared to deduce their similarity, given the first few moments. The calculations were performed using MATLAB for different illustrative values of the parameters. In Figures 2-6 and 8-12 below, the two original PDFs are denoted by solid (red and blue) lines, while the smoothed functions are shown as dashed (red and blue) lines.

\subsection{Case 1: Normal random variables} \label{normal random variables}

For the first illustration, we consider a  normally-distributed random variable $N(\mu,\sigma^2)$ over the compact support $[0,1]$. Its PDF is of the form:

\begin{equation}
    \begin{split}
    f(x) &= \frac{1}{N_f}e^{-\frac{(x-\mu)^2}{2\sigma^2}}\\
    N_f = \int_0^1 e^{-\frac{(x-\mu)^2}{2\sigma^2}} dx &= \sigma \sqrt{\frac{\pi}{2}} \left[ \text{erf}\left( \frac{1-\mu}{\sqrt{2}\sigma} \right) - \text{erf}\left( \frac{-\mu}{\sqrt{2}\sigma} \right) \right]\\
    \end{split}
\end{equation}

\noindent where $N_f$ is the normalization factor, and $\text{erf}(x) = \frac{2}{\sqrt{\pi}}\int_0^x e^{-y^2} dy$ is the error function.

Using \eqref{eq:characteristic function}, we find the characteristic function corresponding to this PDF:

\begin{equation} \label{eq:complicated characteristic function}
    F(k) = e^{-\frac{\mu^2}{2\sigma^2}}e^{-\frac{\sigma^2}{2}\left( k-\frac{i\mu}{\sigma^2} \right)^2} \frac{\left[ \text{erf}\left( \frac{1-(\mu+i\sigma^2 k)}{\sqrt{2}\sigma} \right) - \text{erf}\left( \frac{-(\mu+i\sigma^2 k)}{\sqrt{2}\sigma} \right) \right]}{\left[ \text{erf}\left( \frac{1-\mu}{\sqrt{2}\sigma} \right) - \text{erf}\left( \frac{-x}{\sqrt{2}\sigma} \right) \right]}
\end{equation}

From \eqref{eq:low-pass filter} and \eqref{eq:complicated characteristic function}, and a few variable substitutions, we arrive at an integral expression for the smoothed PDF:

\begin{equation} \label{eq:smoothed function_normal random}
    f^<(x) = f(x) \times \frac{i}{2\sqrt{\pi}} \int_{\frac{\sigma}{\sqrt{2}}\left( \frac{x-\mu}{\sigma^2}+iK \right)}^{\frac{\sigma}{\sqrt{2}}\left( \frac{x-\mu}{\sigma^2}-iK \right)} e^{k^2}\left[ \text{erf}\left( k+\frac{1-x}{\sqrt{2}\sigma} \right) - \text{erf}\left( k+\frac{-x}{\sqrt{2}\sigma} \right) \right] dk
\end{equation}

\begin{figure*}
  \centering
  \begin{tabu} to \textwidth {X[c]}
    \includegraphics[width=.75\linewidth]{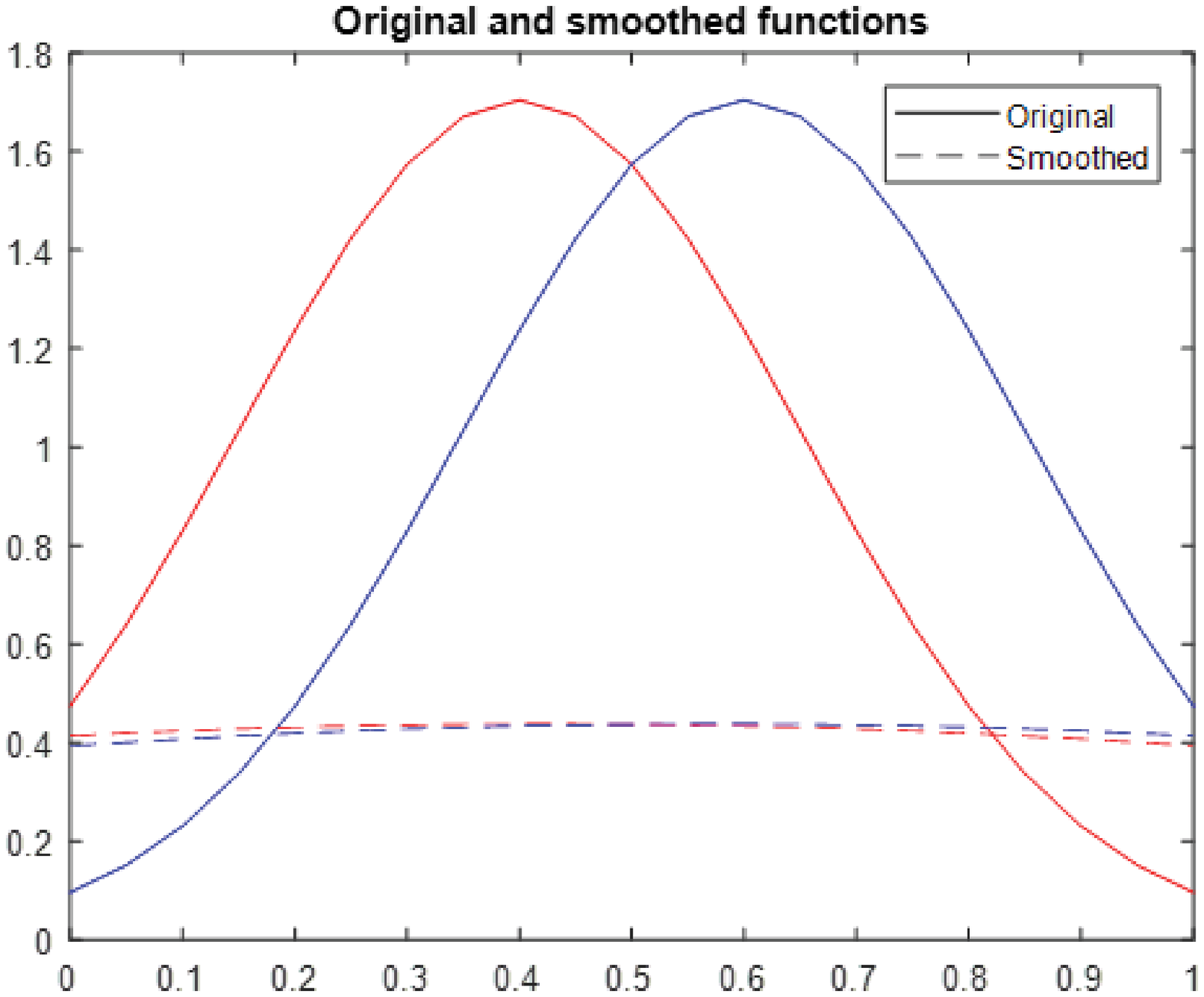} \\
    (a) \\
    \includegraphics[width=.75\linewidth]{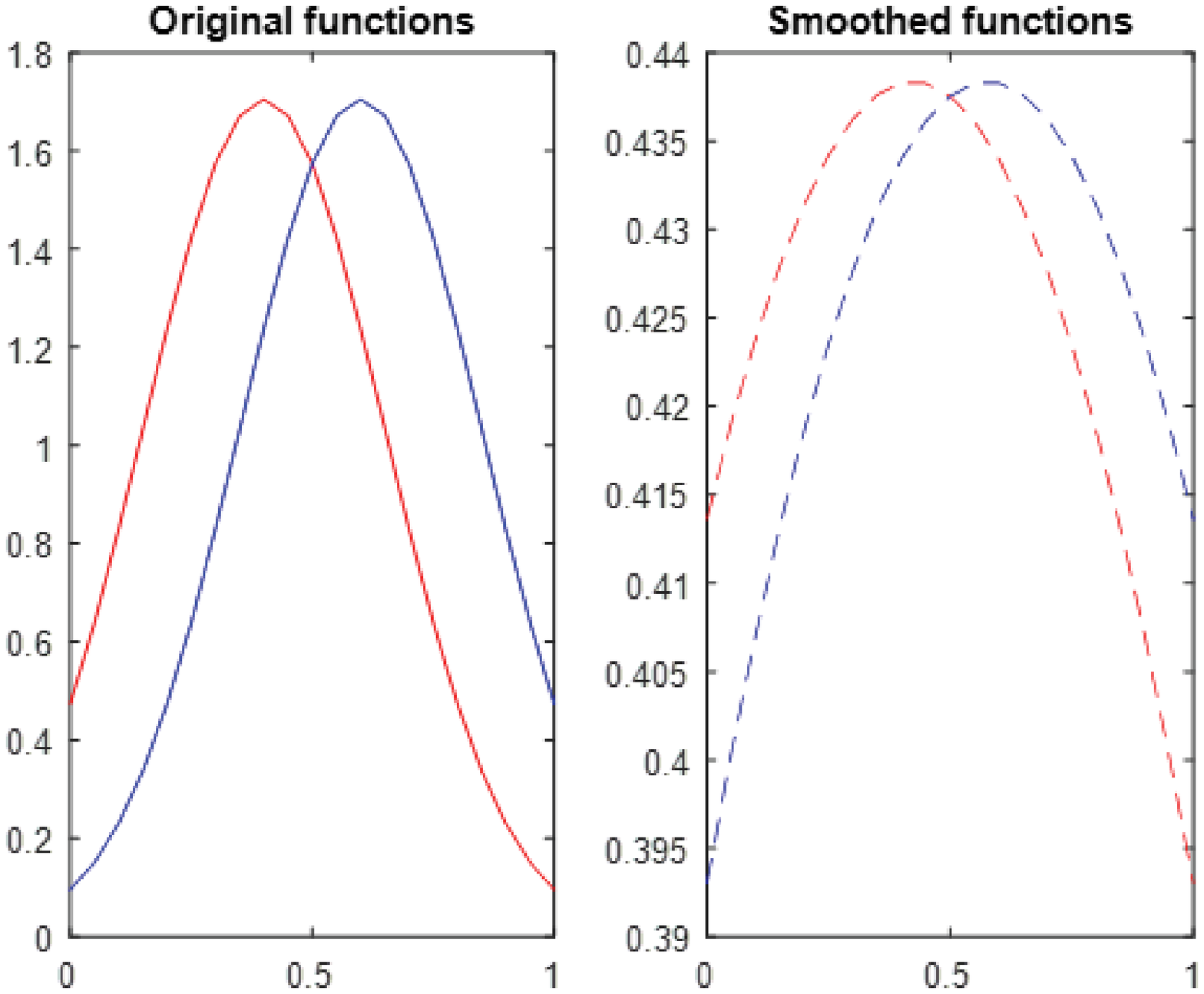} \\
    (b) \\
  \end{tabu}
  \caption{Normally-distributed PDFs supported over $[0,1]$.\\ $n=3$, $\varepsilon=0.1$, Red - $(\mu_1,\sigma_1)=(0.4,0.25)$, Blue - $(\mu_2,\sigma_2)=(0.6,0.25)$\\(a) Comparison of original and smoothed functions\\(b) Smoothed functions shown on different scale for better resolution}
  \label{fig:figure2}
\end{figure*}

\begin{figure*}
  \centering
  \begin{tabu} to \textwidth {X[c]}
    \includegraphics[width=.75\linewidth]{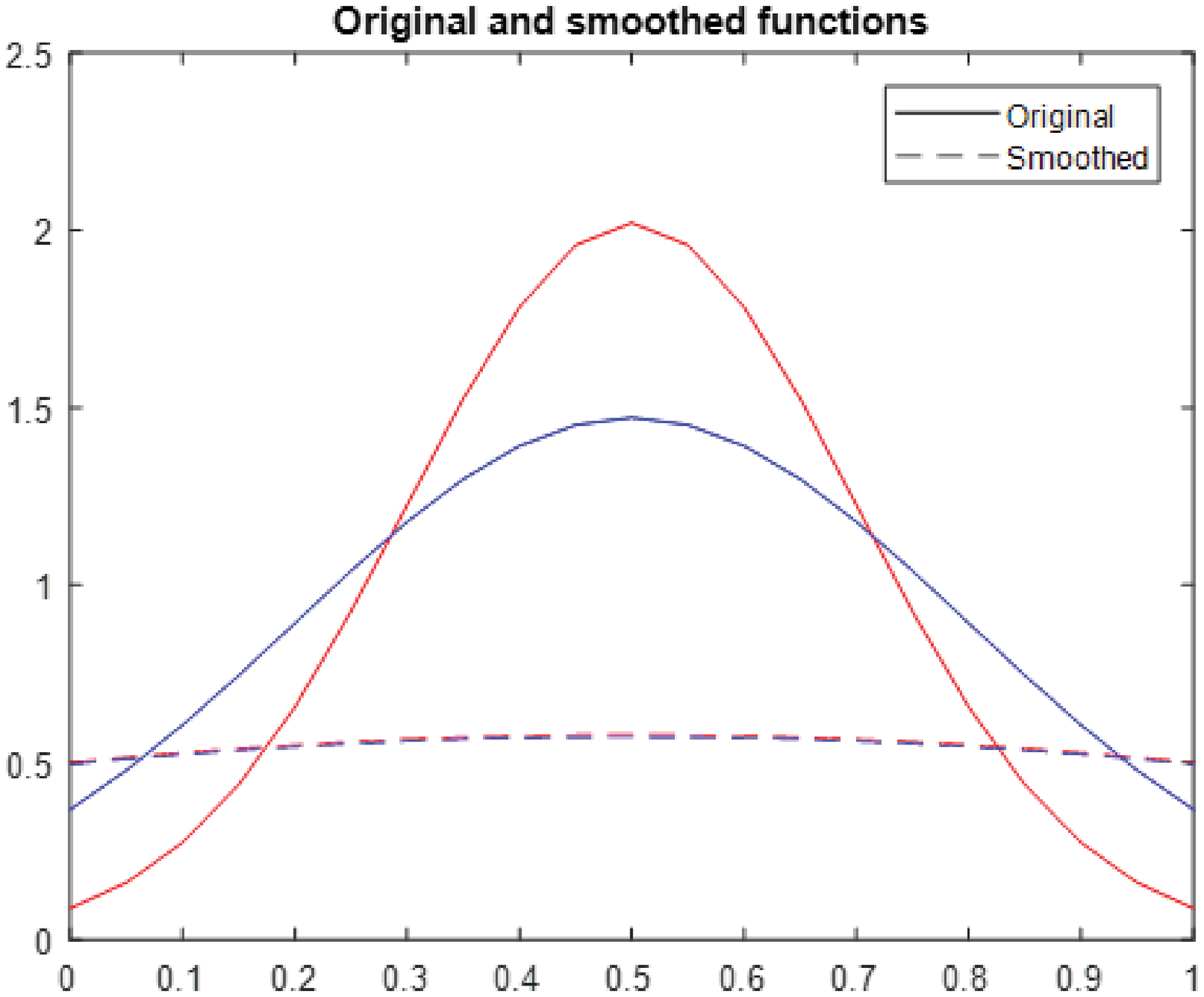} \\
    (a) \\
    \includegraphics[width=.75\linewidth]{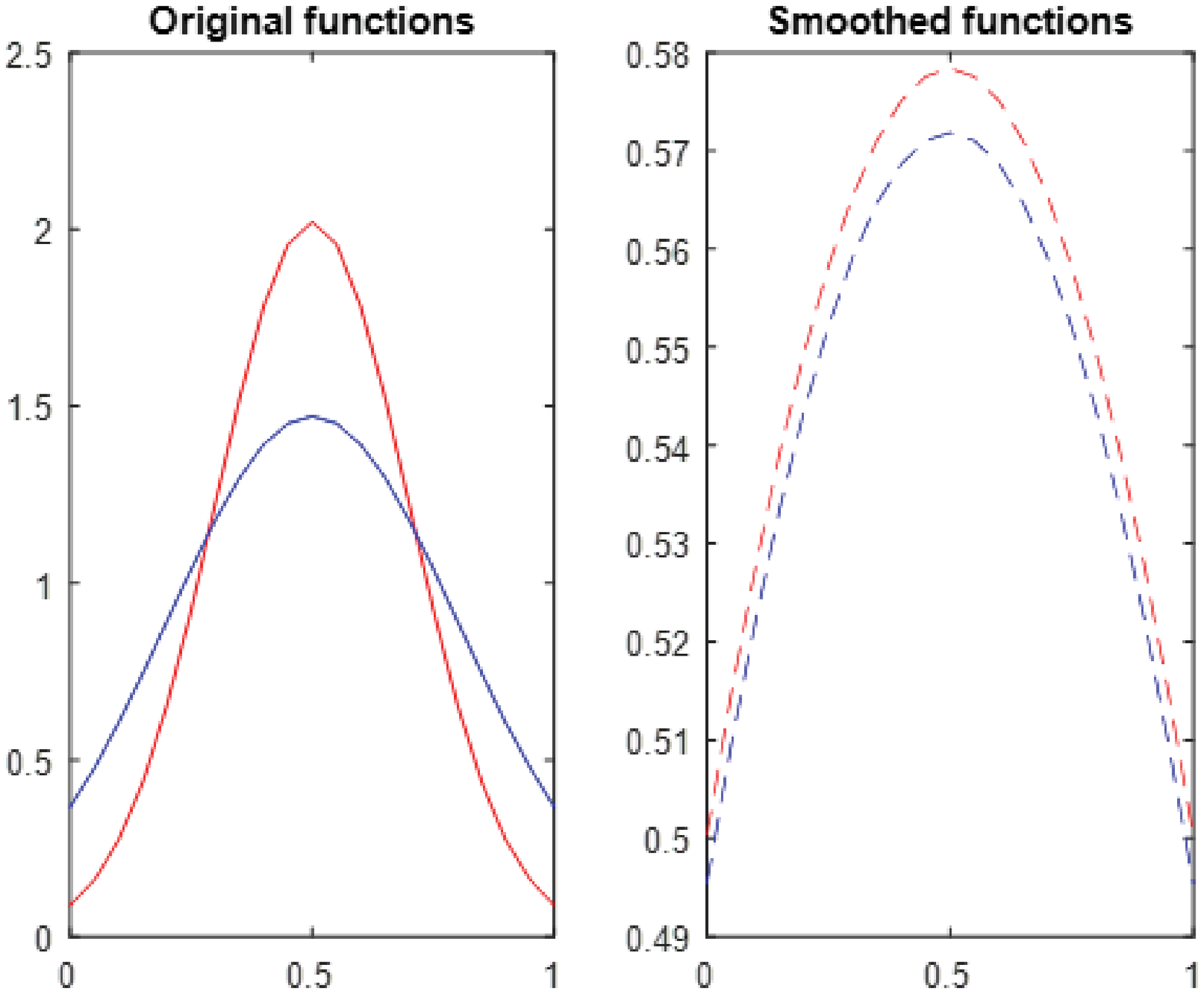} \\
    (b) \\
  \end{tabu}
  \caption[caption]{Normally-distributed PDFs supported over $[0,1]$.\\ $n=3$, $\varepsilon=0.1$, Red - $(\mu_1,\sigma_1)=(0.5,0.2)$, Blue - $(\mu_2,\sigma_2)=(0.5,0.3)$\\(a) Comparison of original and smoothed functions\\(b) Smoothed functions shown on different scale for better resolution}
  \label{fig:figure3}
\end{figure*}

\begin{figure*}
  \centering
  \begin{tabu} to \textwidth {X[c]}
    \includegraphics[width=.75\linewidth]{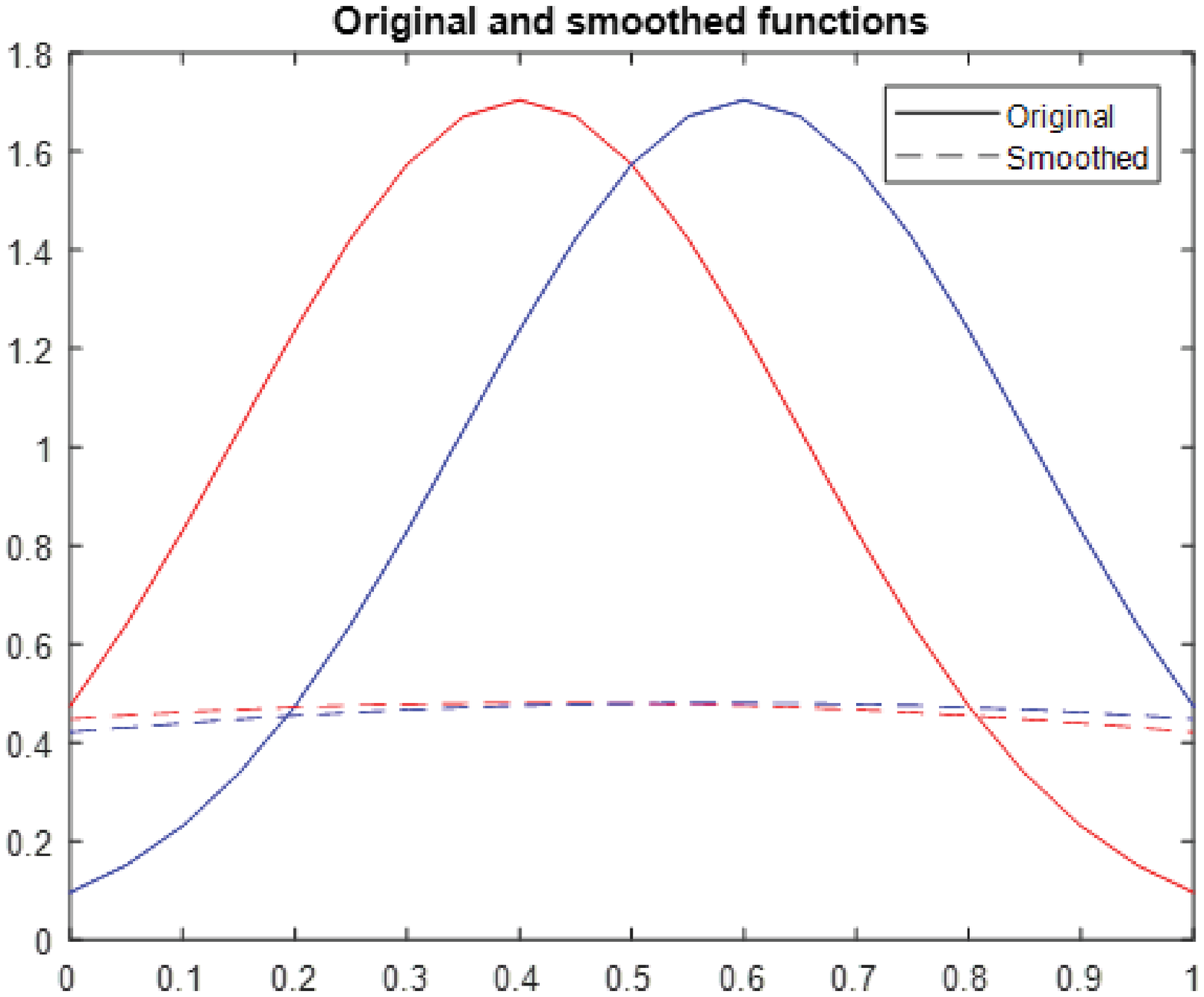} \\
    (a) \\
    \includegraphics[width=.75\linewidth]{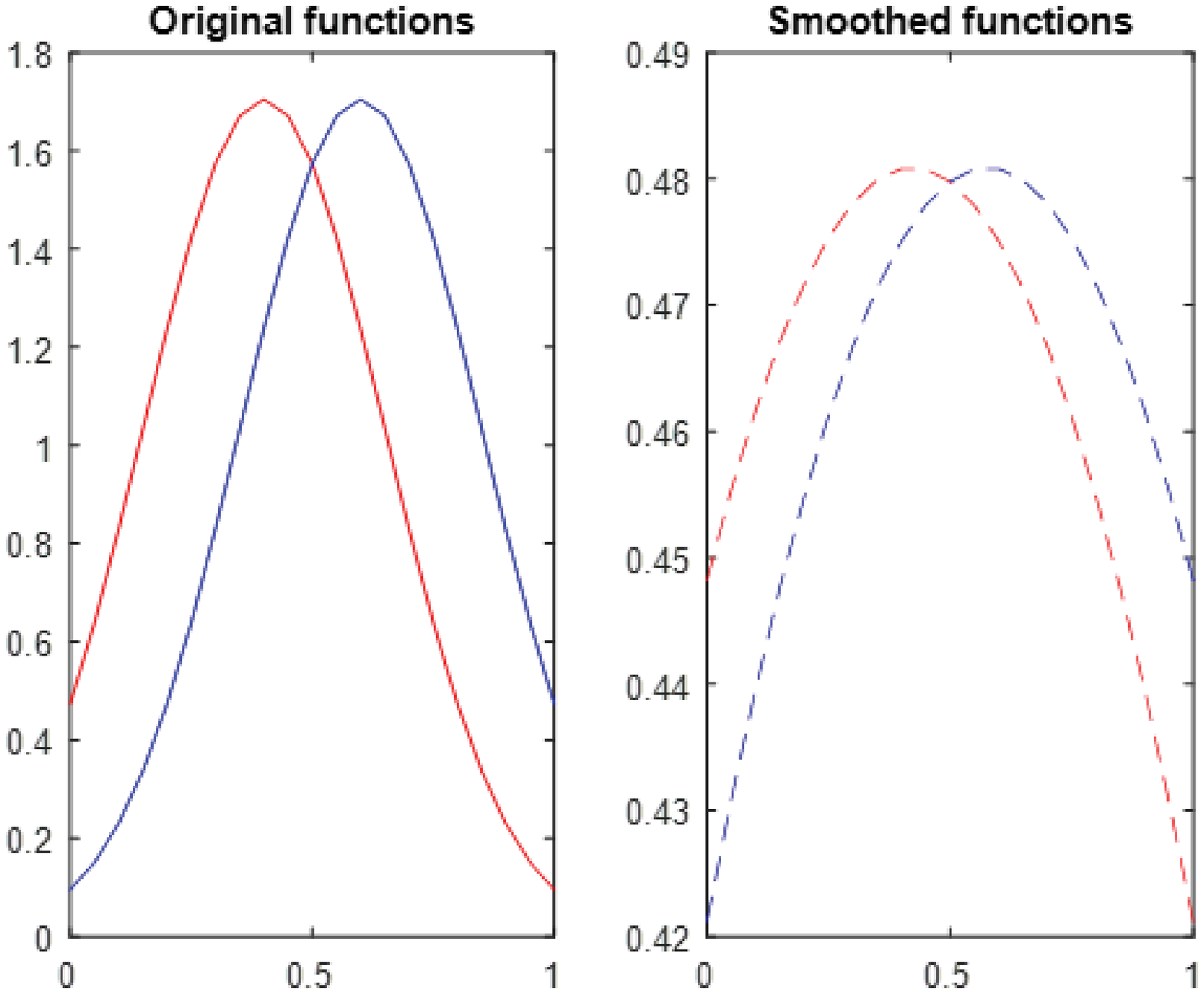} \\
    (b) \\
  \end{tabu}
  \caption{Normally-distributed PDFs supported over $[0,1]$.\\ $n=6$, $\varepsilon=0.1$, Red - $(\mu_1,\sigma_1)=(0.4,0.25)$, Blue - $(\mu_2,\sigma_2)=(0.6,0.25)$\\(a) Comparison of original and smoothed functions\\(b) Smoothed functions shown on different scale for better resolution}
  \label{fig:figure4}
\end{figure*}

\begin{figure*}
  \centering
  \begin{tabu} to \textwidth {X[c]}
    \includegraphics[width=.75\linewidth]{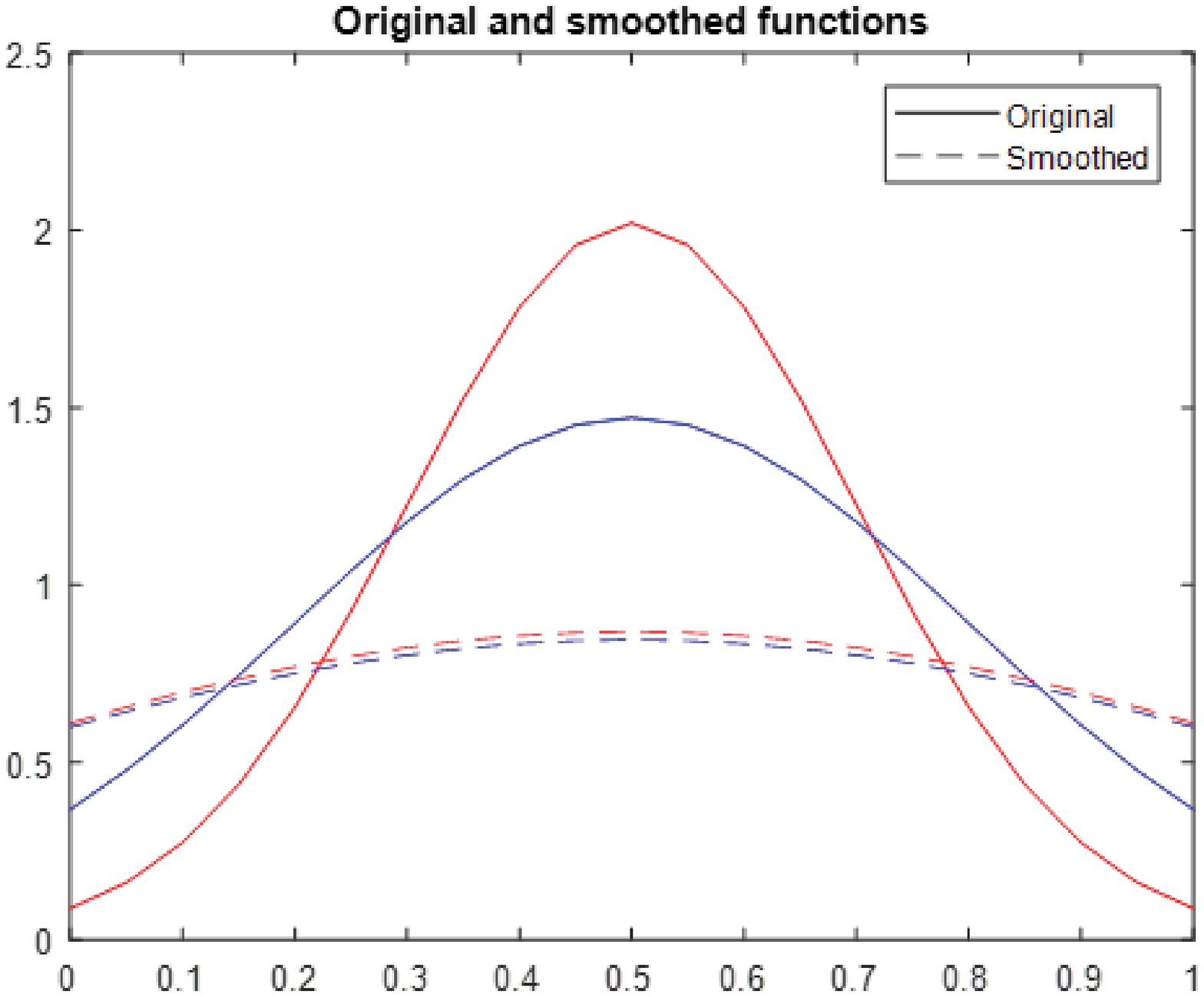} \\
    (a) \\
    \includegraphics[width=.75\linewidth]{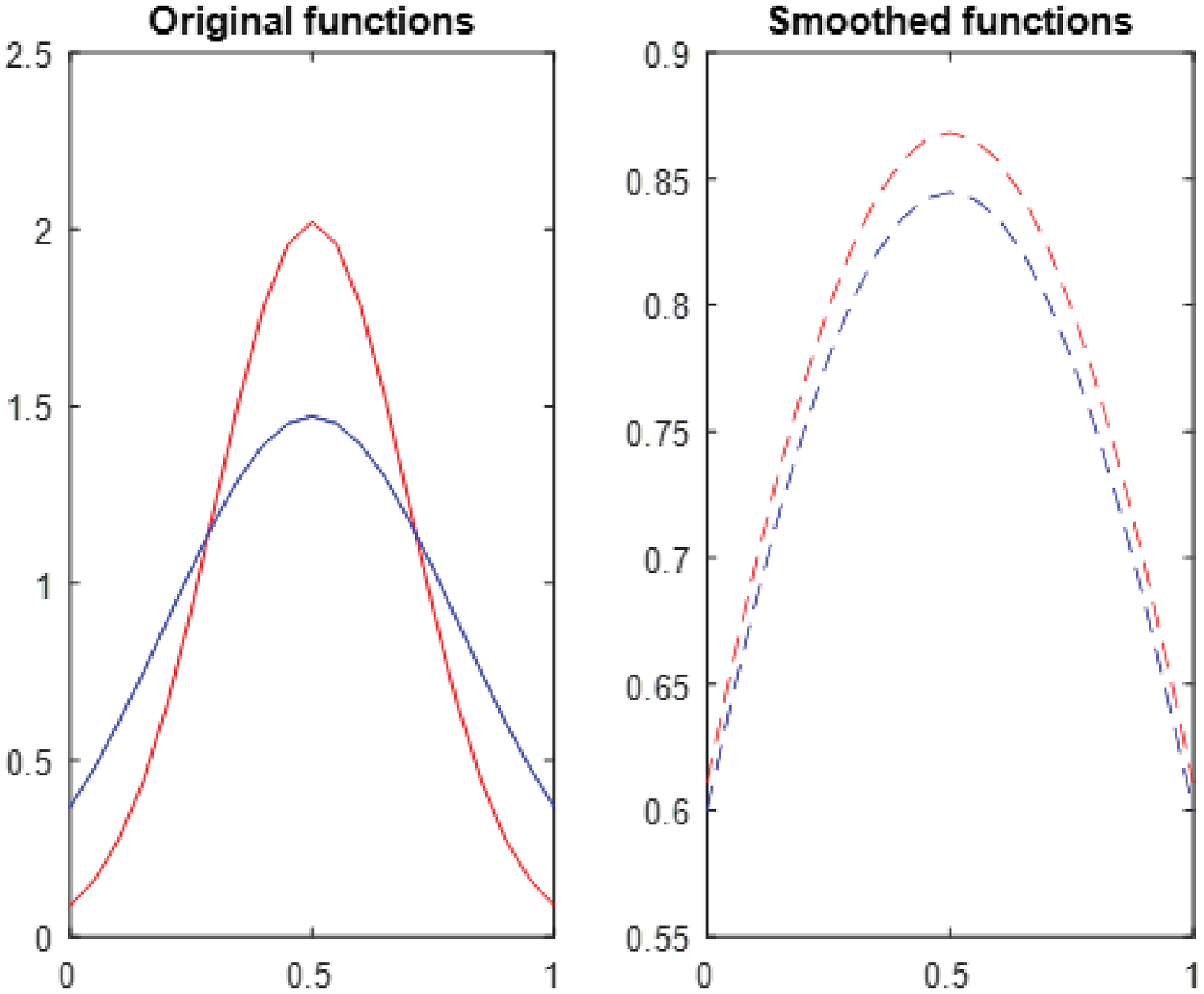} \\
    (b) \\
  \end{tabu}
  \caption[caption]{Normally-distributed PDFs supported over $[0,1]$.\\ $n=6$, $\varepsilon=0.1$, Red - $(\mu_1,\sigma_1)=(0.5,0.2)$, Blue - $(\mu_2,\sigma_2)=(0.5,0.3)$\\(a) Comparison of original and smoothed functions\\(b) Smoothed functions shown on different scale for better resolution}
  \label{fig:figure5}
\end{figure*}

\begin{figure*}
  \centering
  \begin{tabu} to \textwidth {X[c]}
    \includegraphics[width=.75\linewidth]{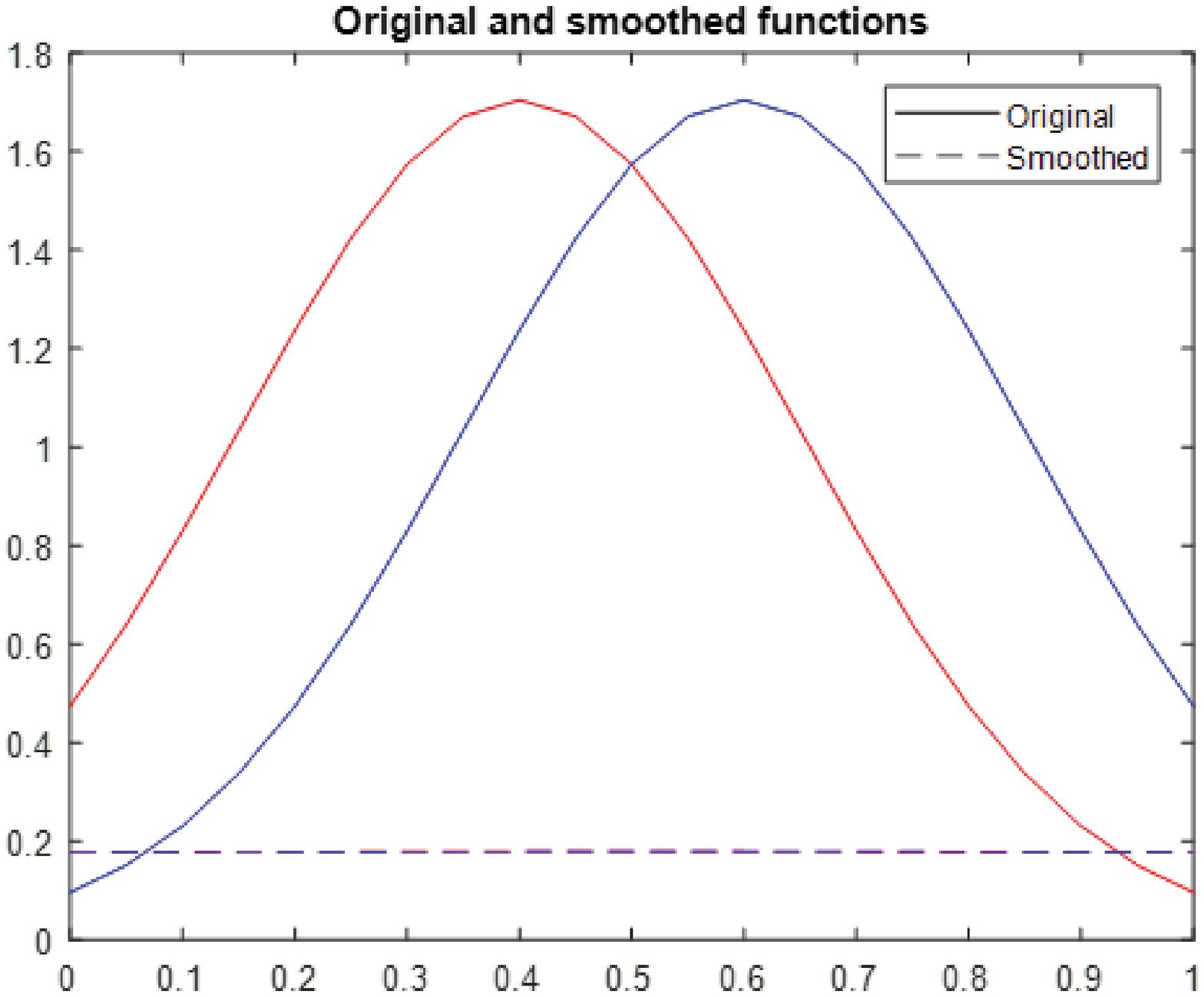} \\
    (a) \\
    \includegraphics[width=.75\linewidth]{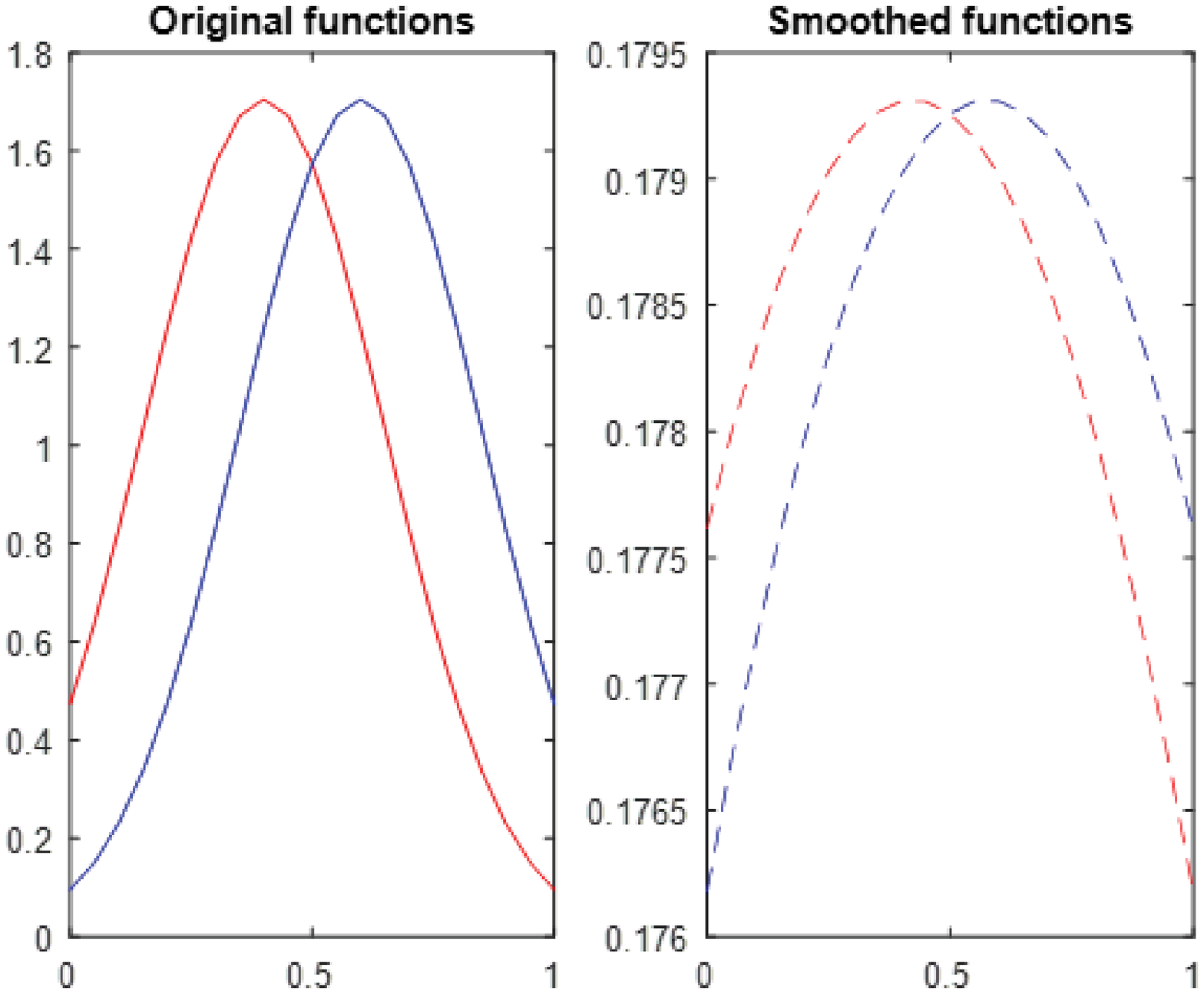} \\
    (b) \\
  \end{tabu}
  \caption{Normally-distributed PDFs supported over $[0,1]$.\\ $n=3$, $\varepsilon=0.01$, Red - $(\mu_1,\sigma_1)=(0.4,0.25)$, Blue - $(\mu_2,\sigma_2)=(0.6,0.25)$\\(a) Comparison of original and smoothed functions\\(b) Smoothed functions shown on different scale for better resolution}
  \label{fig:figure6}
\end{figure*}

\eqref{eq:smoothed function_normal random} was numerically integrated on MATLAB, and the algorithm was applied to two such PDFs and the cut-off scales (and other details) were evaluated. The results are shown in Table \ref{tab:table_normal} and Figures \ref{fig:figure2} - \ref{fig:figure6}. The following observations can be made:

\begin{enumerate}
    \item The cut-off scale is higher (smaller length scales are probed) in the case where $\mu_1=\mu_2,\sigma_1\neq \sigma_2$, when compared to the case of $\mu_1\neq \mu_2,\sigma_1=\sigma_2$. This is expected since it is harder to “match” two PDFs whose peaks are separated, as opposed to when one is simply broader than the other.
    
    \item Increasing the number of moments increases the cut-off (wavenumber) scale. This is because smaller length scales can be probed with more information (moments).

    \item Increasing the number of moments also increases the distance between the smoothed functions. This could perhaps be due to the moments acting as constraints that are to be adhered to while comparing the functions. (It is to be noted that this trend is not observed in the case of the “less-structured” PDFs considered in the next section.)
    
    \item Reducing the tolerance reduces the cut-off scale. Once again, this is expected since a tighter tolerance may be achieved only by sufficient smoothing of the functions.
\end{enumerate}

\subsection{Case 2: Scale-separated PDFs}

Scale-separation is commonly encountered in nature, when the dynamics of a physical system can be separated into two (usually narrow) intervals of length/time scales (see Figure \ref{fig:figure1}). For instance, in certain models of combustion, the assumption of “fast chemistry” is invoked to simplify the analysis. This means that some of the reactions at equilibrium are much faster than others, so that they may be considered instantaneous, and this is a manifestation of temporal scale-separation. As an example of spatial scale-separation, one could consider fluid flow that is not fully turbulent. In such a flow, it is possible to clearly discern the large-scale flow features from the small-scale fluctuations. In summary, scale-separation refers to a case where a range of intermediate (length/time) scales are absent.

\begin{figure}
  \includegraphics[width=\linewidth]{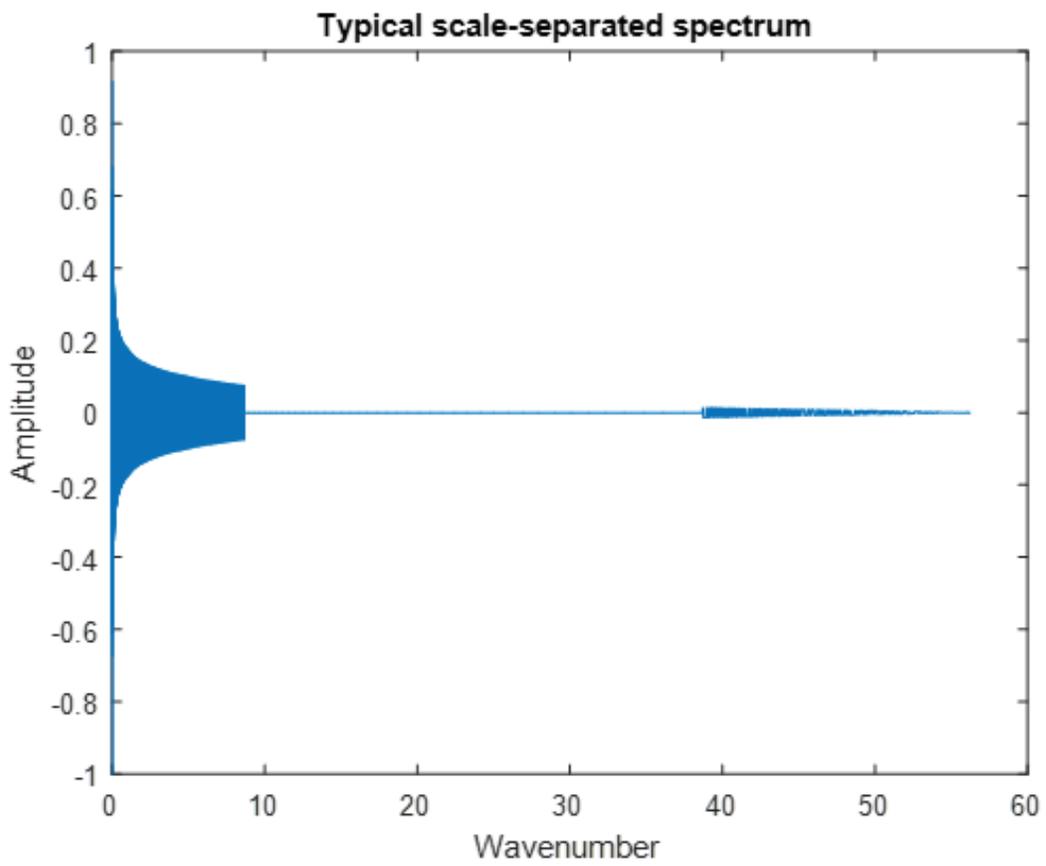}
  \caption{Amplitudes of the sine-component of one of the functions in Figures \ref{fig:figure8} - \ref{fig:figure11}}
  \label{fig:figure7}
\end{figure}

\begin{figure}
  \includegraphics[width=\linewidth]{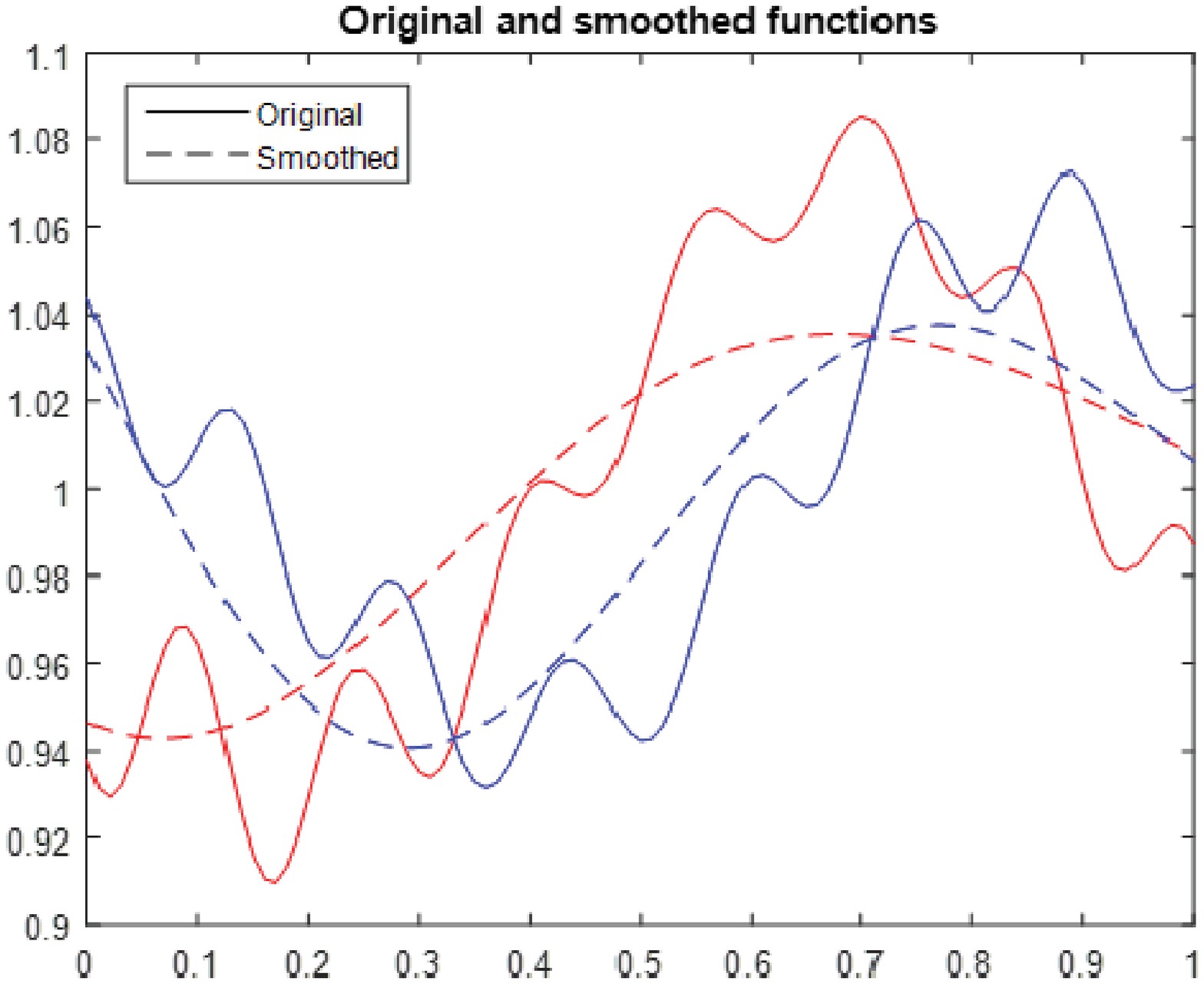}
  \caption{Scale-separated PDFs constructed from random spectra.\\$n=20$, $\varepsilon=1$}
  \label{fig:figure8}
\end{figure}

\begin{figure}
  \includegraphics[width=\linewidth]{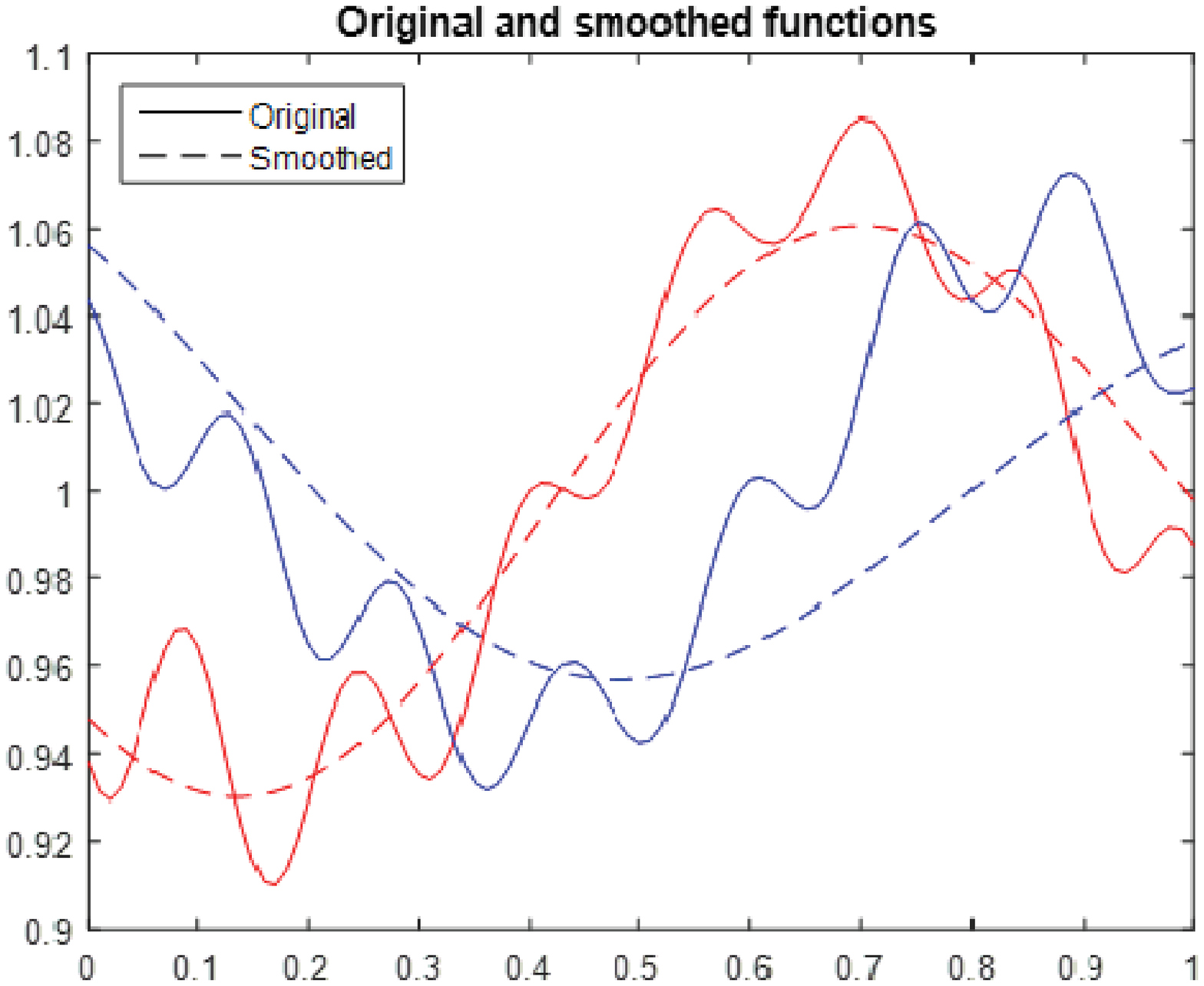}
  \caption{Scale-separated PDFs constructed from random spectra.\\$n=10$, $\varepsilon=1$}
  \label{fig:figure9}
\end{figure}

\begin{figure}
  \includegraphics[width=\linewidth]{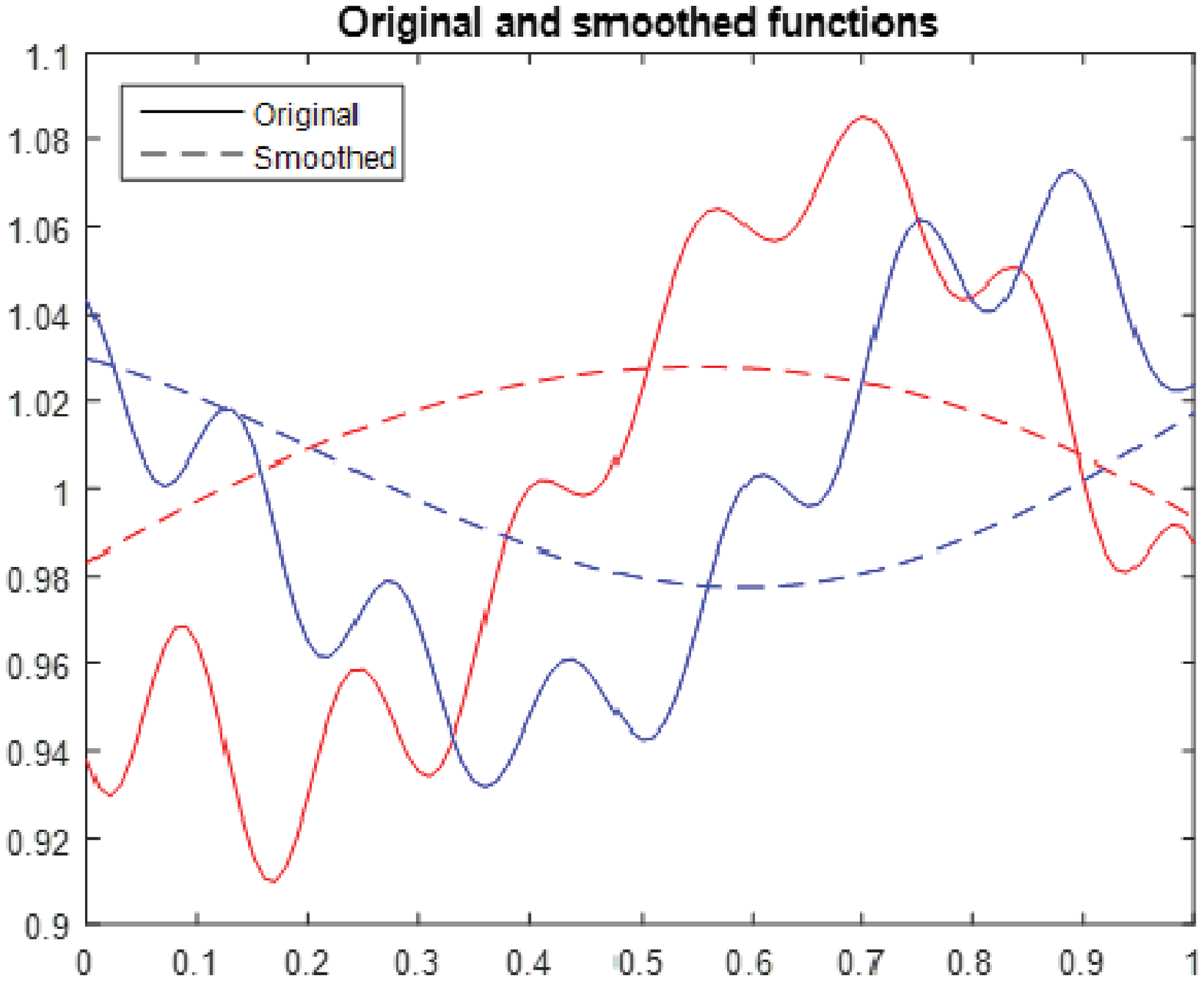}
  \caption{Scale-separated PDFs constructed from random spectra.\\$n=5$, $\varepsilon=1$}
  \label{fig:figure10}
\end{figure}

\begin{figure}
  \includegraphics[width=\linewidth]{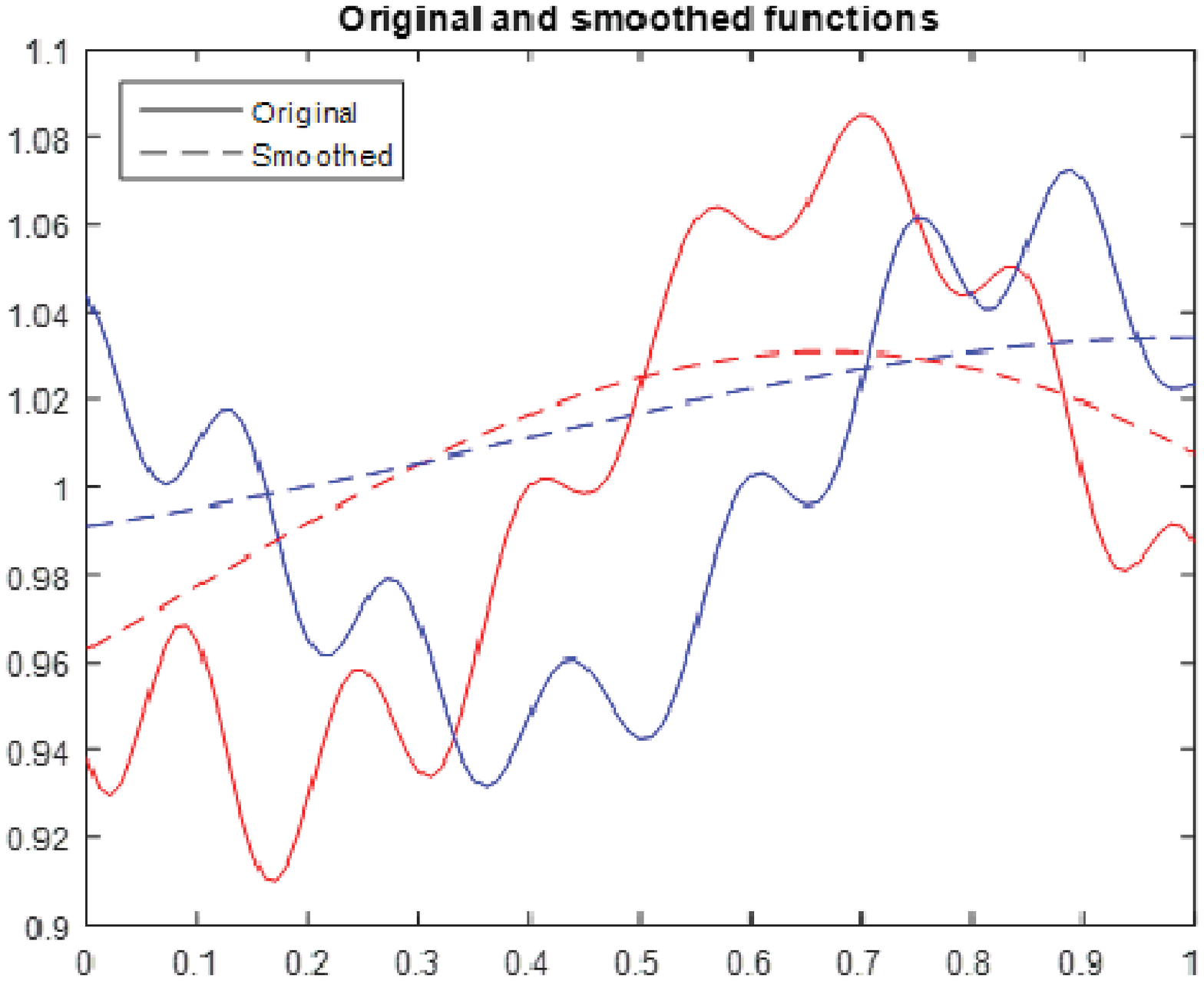}
  \caption{Scale-separated PDFs constructed from random spectra.\\$n=5$, $\varepsilon=0.1$}
  \label{fig:figure11}
\end{figure}

\begin{figure}
  \includegraphics[width=\linewidth]{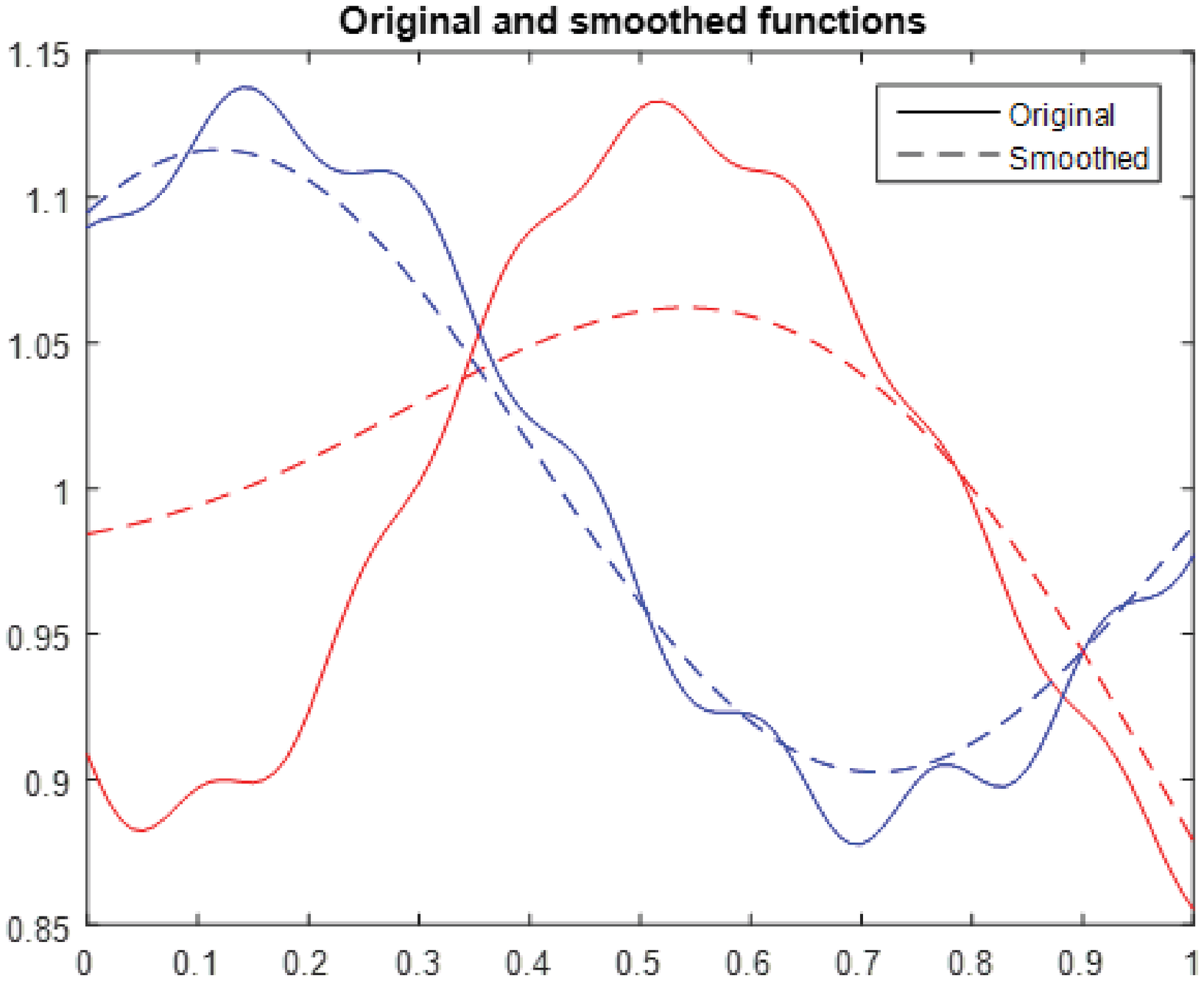}
  \caption{Scale-separated PDFs constructed from random spectra.\\$n=20$, $\varepsilon=1$}
  \label{fig:figure12}
\end{figure}

Due to their ubiquity, it is useful to see an illustration of the method applied to a scale-separated PDF. For this purpose, a spectrum was created by generating normally-distributed random numbers (with zero mean) on MATLAB and arranging them in decreasing order of magnitude. The rearrangement of the numbers was done to imitate the spectrum of a system in a statistically steady state, where fluctuations decrease with decreasing scale size. Further, the spectrum was set to zero in some intermediate wavenumbers to mimic scale-separation (Figure \ref{fig:figure7}). Two such PDFs were constructed as different realizations of the random spectrum and the algorithm described in Section \ref{existence & uniqueness} was used to analyse their similarity and determine the cut-off scale. Using the same random spectrum, the effect of varying parameters is discussed. (The random spectra used to construct Figure \ref{fig:figure12} are different from the ones used for Figures \ref{fig:figure8} - \ref{fig:figure11}.) Integrations were performed numerically and the results are shown in Table \ref{tab:table_scale-separated}.\\

The following observations, similar to those in Section \ref{normal random variables}, can be made:
\begin{enumerate}
    \item Comparing the results for Figures \ref{fig:figure8} and \ref{fig:figure12}, it is seen that the cut-off scale is higher (smaller length scales are probed) in the former. Observing the graphs of the two cases reveals that in Figure \ref{fig:figure12}, the original PDFs are very much “out of phase”, which means they are less similar to begin with than in Figure \ref{fig:figure8}.
    
    \item Increasing the number of moments increases the cut-off (wavenumber) scale.
    
    \item Reducing the tolerance reduces the cut-off scale.
\end{enumerate}

\section{Extension to higher dimensions} \label{extension to higher dimensions}

All of the above analysis was done for PDFs that were functions of just one variable. Can it be extended to a multivariate PDF (i.e., to higher dimensions, say $d$)? It is possible, as will be shown below. The steps are similar to the ones in the 1-dimensional case, but there is a subtle difference in the end.

Let $x = (x_1,x_2,\dots,x_d) \in \mathbb{R}^d$ and let $f(x)$ be a PDF supported over the unit ball in $\mathbb{R}^d$. The characteristic function in $d$ dimensions is given by:

\begin{equation*}
    F(k) = \int_a^b e^{ikx}f(x) dx \ , \ f(x) = \frac{1}{(2\pi)^d}\int_{-\infty}^{\infty} e^{-ikx}F(k) dk
\end{equation*}

\noindent where $k=(k_1,k_2,\dots,k_d) \in \mathbb{R}^d$ is the wavevector. Let $\alpha = (\alpha_1,\alpha_2,\dots,\alpha_d)$ and $|\alpha| = \alpha_1+\alpha_2+\dots+\alpha_d$, where each $\alpha_i \in \{0,1,2,\dots\}$. Then the moments of the PDF are defined as follows:

\begin{equation}
    \langle x^{\alpha} \rangle := \langle x_1^{\alpha_1}x_2^{\alpha_2}\dots x_d^{\alpha_d}\rangle = \int x_1^{\alpha_1}x_2^{\alpha_2}\dots x_d^{\alpha_d}f(x) dx = (-i)^{|\alpha|}D^{\alpha}F(0)
\end{equation}

\noindent where $D^{\alpha} \equiv \frac{\partial^{\alpha}}{\partial k_1^{\alpha 1} \partial k_2^{\alpha 2} \dots \partial k_d^{\alpha d}}$.

The upper bounds for the moments are defined as:

$$M_{\alpha} \ge |\langle x^{\alpha} \rangle_f - \langle x^{\alpha} \rangle_g| = |D^{\alpha}F(0) - D^{\alpha}G(0)|$$

The low-pass filtering operator is defined for a $d$-dimensional cut-off scale $K = (K_1,K_2,\dots,K_d)$:

\begin{equation}
    P_K[f(x)] = f^<(x) = \frac{1}{(2\pi)^d}\int_{-K_d}^{K_d}\dots \int_{-K_2}^{K_2}\int_{-K_1}^{K_1} e^{-ikx}F(k) dk_1dk_2\dots dk_d
\end{equation}

Expanding the characteristic function in a Taylor series about the origin,

\begin{equation}
    F(k) = \sum_{|\alpha|\le n-1} \frac{D^{\alpha}F(0)}{\alpha !}k^{\alpha} + \sum_{|\alpha|= n} \frac{D^{\alpha}F(tk)}{\alpha !}k^{\alpha}
\end{equation}

\noindent where $t\in (0,1)$ and $\alpha ! = \alpha_1 ! \alpha_2 !\dots \alpha_d !$.

As in the 1-dimensional case, the remainder term in the Taylor series can be bounded as:

$$|D^{\alpha}F(tk)| \le \langle x^{\alpha} \rangle \le \langle x^{\beta} \rangle$$

$\beta$ is obtained by reducing any one non-zero component of $\alpha$ by $1$. Finally, the distance metric in $d$-dimension is:

\begin{equation}
    d(f,g) = \frac{\int_{B(0,1)}|f(x)-g(x)|dx}{\int_{B(0,1)}dx}
\end{equation}

Following the same sequence of inequalities as before, the distance inequality for the smoothed functions becomes:

\begin{equation*}
    d(f^<,g^<) \le \frac{1}{(2\pi)^d} \int_{-K_d}^{K_d}\dots \int_{-K_2}^{K_2}\int_{-K_1}^{K_1} \left\{ \sum_{|\alpha| \le n-1} \frac{M_{\alpha}}{\alpha!}|k|^{\alpha} + \sum_{|\alpha| = n} \frac{R_{\alpha}}{\alpha!}|k|^{\alpha} \right\} dk_1dk_2\dots dk_d
\end{equation*}

\noindent where $R_{\alpha}=\langle x^{\beta_f} \rangle_f + \langle x^{\beta_g} \rangle_g$ is the remainder term. Here, $\beta_f$ and $\beta_g$ are chosen so that the remainder term is the least possible (among all the moments of that order). Note that $|k|^{\alpha} = |k_1|^{\alpha_1}|k_2|^{\alpha_2}\dots |k_d|^{\alpha_d}$.

\begin{equation}
    \therefore d(f^<,g^<) \le \frac{1}{\pi^d} \left\{ \sum_{|\alpha|\le n-1} \frac{M_{\alpha}}{(\alpha + \mathbb{I})!}K^{\alpha + \mathbb{I}} + \sum_{|\alpha|=n} \frac{R_{\alpha}}{(\alpha + \mathbb{I})!}K^{\alpha + \mathbb{I}} \right\} = \varepsilon
\end{equation}

\noindent where $\mathbb{I} = (1,1,\dots,1)$.

The difference in the higher dimensional case $(d>1)$ is that the solution $(K_1,K_2,…,K_d)$ is not unique. However, a unique (and conservative) value for the cut-off scale, that is common over all the dimensions, can be determined. If $\max_{i} \{K_i \}=\kappa$, then setting all the $K_i=\kappa$, the resulting equation is of the form:

\begin{equation}
    \frac{\kappa^d}{\pi^d} \left\{ \sum_{j=0}^{n} a_j \kappa_j \right\} = \varepsilon
\end{equation}

\noindent where the $a_j$ are the sum of various moments (divided by the appropriate factorials) and the remainder term. The existence and uniqueness of this $\kappa$ can be proved just as in the 1-dimensional case.

\section{Conclusions} \label{conclusions}

A method has been proposed to estimate a length scale to which given PDFs must be smoothed for their (normalized) $L^1$ distance to be less than a user-specified tolerance, given the first few moments of the PDFs. It has been shown that such a cut-off scale indeed exists and is unique. Two numerical examples were used as proof of concept, as well as to illustrate the working of the algorithm. An extension to higher dimensions has been outlined. This scheme is hoped to find use in data analysis for comparing large data sets, as calculating moments is less computationally-intensive than having to reverse-engineer the PDFs of the data sets in order to compare them.

\section{Scope for future work} \label{scope for future work}

The calculations in Section \ref{estimates} involves a series of inequalities, which lead to a very conservative estimate of the distance between the smoothed functions and the cut-off scale, as seen in Tables \ref{tab:table_normal} and \ref{tab:table_scale-separated}. A useful direction for future research is to sharpen the estimates with more accurate inequalities.

Moreover, the method described in this article only deals with PDFs having bounded support. As discussed in Section \ref{unbounded support}, unbounded domains pose major problems in the existence or the finiteness of moments. Extending this scheme, or formulating an entirely new one, for unbounded domains could be another meaningful and interesting research problem.

\ack PCJ wishes to thank Andy Sebastian and Andrew Corson for helpful discussions in the preliminary stages of this research and acknowledge Joydeep Singha's suggestion to extend this algorithm to multiple dimensions. The authors are also grateful to Dr. Chris Jarzynski and Dr. Venkatarathnam Gadhiraju for their valuable suggestions to improve the manuscript.

\appendix

\section{Bounds on higher moments} \label{bounds on higher moments}

Consider a PDF $f(x)$ with compact support $[0,1]$. Then, its $n^{\text{th}}$ moment can be bounded from above by a lower moment $m$ and/or other terms, depending on the smoothness of the PDF. This is useful in bounding the remainder term in the algorithm discussed above.

\subsection{PDF is bounded: $f(x) \in L^{\infty}([0,1])$}

For all $m,n \in \mathbb{N}$ with $m<n$, we use the Cauchy-Schwartz inequality to obtain:

\begin{align*}
    \langle x^n \rangle = \int_0^1 x^n f(x) dx &\le \left( \int_0^1 x^{2(n-m)} dx \right)^{\frac{1}{2}} \left( \int_0^1 x^{2m} f^2(x) dx \right)^{\frac{1}{2}}\\
    &\le \left( \frac{1}{2(n-m)+1} \right)^{\frac{1}{2}} \left( \left\Vert f \right\Vert_{L^{\infty}} \int_0^1 x^{2m} f(x) dx \right)^{\frac{1}{2}}
\end{align*}

\begin{equation} \label{appendix_bounded}
    \langle x^n \rangle \le \left( \frac{\left\Vert f \right\Vert_{L^{\infty}} \langle x^{2m} \rangle}{2(n-m)+1} \right)^{\frac{1}{2}}
\end{equation}

\noindent Three special cases of \eqref{appendix_bounded} may be of interest.

\begin{enumerate}
    \item $m=0$
    
    $$\Rightarrow \langle x^n \rangle \le \left( \frac{\left\Vert f \right\Vert_{L^{\infty}}}{2n+1} \right)^{\frac{1}{2}}$$
    
    \item $n$ is even and $m=\frac{n}{2}$
    
    $$\Rightarrow \langle x^n \rangle \le \frac{\left\Vert f \right\Vert_{L^{\infty}}}{n+1}$$
    
    \item $n$ is odd and $m=\frac{n-1}{2}$
    
    $$\Rightarrow \langle x^n \rangle \le \left( \frac{\left\Vert f \right\Vert_{L^{\infty}} \langle x^{n-1} \rangle}{n+2} \right)^{\frac{1}{2}}$$
\end{enumerate}

\subsection{PDF is absolutely continuous}

In this case, for some $a \in [0,1]$,

$$|f(x)| \le |f(a)| + \left\vert \int_a^x |f'(y)|dy \right\vert \le |f(a)| + \int_0^1 |f'(y)|dy$$

\begin{equation} \label{appendix_absolutely continuous}
    \Rightarrow \left\Vert f \right\Vert_{L^{\infty}} \le |f(a)| + \left\Vert f' \right\Vert_{L^1} \le |f(a)| + \left\Vert f' \right\Vert_{L^p}
\end{equation}

\noindent for $1\le p\le \infty$. This bound may be used in \eqref{appendix_bounded}.

\subsection{$f(1), f^{(1)}(1), \left\Vert f^{(2)} \right\Vert_{L^{\infty}}$ are finite and known}

\begin{align*}
    \langle x^n \rangle &= \int_0^1 x^n f(x) dx  \xrightarrow[\text{parts, twice}]{\text{integrating by}} \\
    &= \frac{f(1)}{n+1} - \frac{f^{(1)}(1)}{(n+1)(n+2)} + \frac{\int_0^1 x^{n+2}f^{(2)}(x) dx}{(n+1)(n+2)}\\
\end{align*}

\begin{equation}
    \Rightarrow \langle x^n \rangle \le \frac{f(1)}{n+1} - \frac{f^{(1)}(1)}{(n+1)(n+2)} + \frac{\left\Vert f^{(2)} \right\Vert_{L^{\infty}}}{(n+1)(n+2)(n+3)}
\end{equation}

This can, of course, be extended to higher derivatives by integrating by parts repeatedly, depending on how much information one already has about the PDF.

\section{Smoothness of the characteristic function} \label{appendix_smoothness_characteristic function}

\begin{thm}
    If a random variable has moments up to the order $n$, then the corresponding characteristic function belongs to $C^n(\mathbb{R})$.
\end{thm}

\textit{Proof}: (Adapted from \cite{Lukacs1970CharacteristicFunctions}) From \eqref{eq:characteristic function}, we know that for a random variable with a PDF given by $f(x)$,

\begin{equation*}
    F^{(n)}(k) = i^n \int_0^1 x^n f(x) e^{ikx} dx 
\end{equation*}

In the case that a PDF does not exist, this can be rewritten in terms of the cumulative distribution function as:

\begin{equation*}
    F^{(n)}(k) = i^n \int_0^1 x^n e^{ikx} d\Phi(x) 
\end{equation*}

Thus, if the PDF exists, then $f(x) dx = d\Phi(x)$. Now, consider 

\begin{equation*}
    \left\vert F^{(n)}(k+h) - F^{(n)}(k) \right\vert \le \int_0^1 x^n \left\vert e^{ihx} -1 \right\vert d\Phi(x) = 2\int_0^1 x^n \left\vert \sin{\frac{hx}{2}} \right\vert d\Phi(x)
\end{equation*}

The RHS is independent of $k$, and less than $2\int_0^1 x^n d\Phi(x)$ which is twice the $n^{\text{th}}$ moment. Also, the RHS can be made arbitrarily small by taking the limit $h \rightarrow 0$. From all these observations, we conclude that $F^{(n)}(k)$ is uniformly continuous on the entire real line. Thus, $F(k) \in C^n(\mathbb{R})$.

\qed \\

From the discussion in Section \ref{bounded support}, we know that a random variable over a compact support has all moments. Combining this with the above theorem gives us the following corollary.

\begin{cor}
    For a random variable over a compact support, the characteristic function belongs to $C^{\infty}(\mathbb{R})$.
\end{cor}

\section*{References}

\bibliographystyle{unsrt}
\bibliography{mendeley_v2}

\noappendix

\Tables

\begin{center}
\begin{table}
\caption{\label{label}Effect of various parameters for a normally-distributed PDF over a compact support.}
\begin{center}
\begin{tabu} to 1.0\textwidth {|X[c,m]|X[c,m]|X[c,m]|X[c,m]|X[c,m]|X[c,m]|X[c,m]|}
    \hline
    Figure number & Number of moments $(n)$ & Tolerance $(\varepsilon)$ & Cut-off scale $(K)$ & Original distance $d(f,g)$ & Smoothed distance $d(f^<,g^<)$ & Parameters of PDFs $(\mu_1, \sigma_1)$ $(\mu_2, \sigma_2)$\\
    \hline\hline
 
    \hline
    Figure \ref{fig:figure2} & 3 & 0.1 & 1.396927 & 0.561495	& 0.010545 & $(0.4,0.25)$ $(0.6,0.25)$\\
    \hline
    
    \hline
    Figure \ref{fig:figure3} & 3 & 0.1 & 1.854235 & 0.289393	& 0.005971 & $(0.5,0.2)$ $(0.5,0.3)$\\
    \hline
    
    \hline
    Figure \ref{fig:figure4} & 6 & 0.1 & 1.537326 & 0.561495	& 0.013916 & $(0.4,0.25)$ $(0.6,0.25)$\\
    \hline
    
    \hline
    Figure \ref{fig:figure5} & 6 & 0.1 & 2.863850 & 0.289393	& 0.019028 & $(0.5,0.2)$ $(0.5,0.3)$\\
    \hline
    
    \hline
    Figure \ref{fig:figure6} & 3 & 0.01 & 0.564254 & 0.561495	& 0.000723 & $(0.4,0.25)$ $(0.6,0.25)$\\
    \hline
    
\end{tabu}
\end{center}
\label{tab:table_normal}
\end{table}
\end{center}

\begin{table}
\centering
\caption{\label{label}Effect of various parameters for a scale-separated PDF over a compact support.}
\begin{center}
\begin{tabu} to \textwidth {|X[c,m]|X[c,m]|X[c,m]|X[c,m]|X[c,m]|X[c,m]|}
    \hline
    Figure number & Number of moments $(n)$ & Tolerance $(\varepsilon)$ & Cut-off scale $(K)$ & Original distance $d(f,g)$ & Smoothed distance $d(f^<,g^<)$\\
    \hline\hline
 
    \hline
    Figure \ref{fig:figure8} & 20 & 1 & 6.949457 & 0.048280	& 0.024368 \\
    \hline
    
    \hline
    Figure \ref{fig:figure9} & 10 & 1 & 6.042744 & 0.048280	& 0.057852 \\
    \hline
    
    \hline
    Figure \ref{fig:figure10} & 5 & 1 & 4.178820 & 0.048280	& 0.029120 \\
    \hline
    
    \hline
    Figure \ref{fig:figure11} & 5 & 0.1 & 2.816885 & 0.048280	& 0.009402 \\
    \hline
    
    \hline
    Figure \ref{fig:figure12} & 20 & 1 & 6.364578 & 0.135750	& 0.086696 \\
    \hline
    
\end{tabu}
\end{center}
\label{tab:table_scale-separated}
\end{table}

\noappendix

\end{document}